\documentclass[10pt]{amsart}
\usepackage[cp1251]{inputenc}
\usepackage[english]{babel}
\usepackage{amsmath}
\usepackage{amssymb}
\usepackage{amsfonts}

\begin{document}
\renewcommand{\refname}{References}

\thispagestyle{empty}

\title[Explicit One-Step Strong Numerical Methods of Orders 2.0 and 2.5]
{Explicit One-Step Strong Numerical Methods of Orders 2.0 and 2.5 for
Ito Stochastic Differential Equations Based on the 
Unified Taylor--Ito and Taylor--Stratonovich
Expansions}
\author[D.F. Kuznetsov]{Dmitriy F. Kuznetsov}
\address{Dmitriy Feliksovich Kuznetsov
\newline\hphantom{iii} Peter the Great Saint-Petersburg Polytechnic University,
\newline\hphantom{iii} Polytechnicheskaya ul., 29,
\newline\hphantom{iii} 195251, Saint-Petersburg, Russia}%
\email{sde\_kuznetsov@inbox.ru}

\maketitle {\small
\begin{quote}
\noindent{\sc Abstract.} 
The article is devoted to the construction of 
explicit one-step strong numerical methods
with the orders 2.0 and 2.5 of convergence for Ito
stochastic differential equations with
multidimensional non-commutative noise.
We consider numerical methods based
on the unified Taylor--Ito and Taylor--Stratonovich
expansions.
For the numerical modeling of iterated Ito and Stratonovich
stochastic integrals of multiplicities 1 to 5 we apply the method
of multiple Fourier--Legendre series converging in the
sense of norm in Hilbert space $L_2([t, T]^k),$ $k=1,\ldots,5$.
The article is addressed to engineers who use
numerical modeling in stochastic control
and for solving the non-linear filtering problem. 
The article will be interesting 
to scientists who working in the field 
of numerical integration of stochastic differential equations.
\medskip
\end{quote}
}

\vspace{12mm}


\setlength{\baselineskip}{1.8em}

\tableofcontents

\setlength{\baselineskip}{1.2em}


\section{Introduction}

\vspace{5mm}

Let $(\Omega,$ ${\rm F},$ ${\sf P})$ be a complete probability space, let 
$\{{\rm F}_t, t\in[0,T]\}$ be a nondecreasing 
right-continous family of $\sigma$-subfields of ${\rm F},$
and let ${\bf f}_t$ be a standard $m$-dimensional Wiener 
stochastic process, which is
${\rm F}_t$-measurable for any $t\in[0, T].$ We assume that the components
${\bf f}_{t}^{(i)}$ $(i=1,\ldots,m)$ of this process are independent. Consider
an Ito stochastic differential equation in the integral form

\vspace{-1mm}
\begin{equation}
\label{1.5.2}
{\bf x}_t={\bf x}_0+\int\limits_0^t {\bf a}({\bf x}_{\tau},\tau)d\tau+
\int\limits_0^t B({\bf x}_{\tau},\tau)d{\bf f}_{\tau},\ \ \
{\bf x}_0={\bf x}(0,\omega).
\end{equation}

\vspace{2mm}
\noindent
Here ${\bf x}_t$ is some $n$-dimensional stochastic process 
satisfying the equation (\ref{1.5.2}). 
The nonrandom functions ${\bf a}: \mathbb{R}^n\times[0, T]\to\mathbb{R}^n$,
$B: \mathbb{R}^n\times[0, T]\to\mathbb{R}^{n\times m}$
guarantee the existence and uniqueness up to stochastic equivalence 
of a solution
of the equation (\ref{1.5.2}) \cite{1}. The second integral on the 
right-hand side of (\ref{1.5.2}) is 
interpreted as the Ito stochastic integral.
Let ${\bf x}_0$ be an $n$-dimensional random variable, which is 
${\rm F}_0$-measurable and 
${\sf M}\{\left|{\bf x}_0\right|^2\}<\infty$ 
(${\sf M}$ denotes a mathematical expectation).
We assume that
${\bf x}_0$ and ${\bf f}_t-{\bf f}_0$ are independent when $t>0.$

It is well known \cite{KlPl2}-\cite{Mi3}
that Ito stochastic differential equations are 
adequate mathematical models of dynamic systems under 
the influence of random disturbances. One of the effective approaches 
to numerical integration of 
Ito stochastic differential equations is an approach based on 
the Taylor--Ito and 
Taylor--Stratonovich expansions
\cite{KlPl2}-\cite{arxiv-24999}. The most important feature of such 
expansions is a presence in them of the so-called iterated
Ito and Stratonovich stochastic integrals, which play the key 
role for solving the 
problem of numerical integration of Ito stochastic 
differential equations and have the 
following form

\vspace{-1mm}
\begin{equation}
\label{ito}
J[\psi^{(k)}]_{T,t}=\int\limits_t^T\psi_k(t_k) \ldots \int\limits_t^{t_{2}}
\psi_1(t_1) d{\bf w}_{t_1}^{(i_1)}\ldots
d{\bf w}_{t_k}^{(i_k)},
\end{equation}

\begin{equation}
\label{str}
J^{*}[\psi^{(k)}]_{T,t}=
{\int\limits_t^{*}}^T
\psi_k(t_k) \ldots 
{\int\limits_t^{*}}^{t_2}
\psi_1(t_1) d{\bf w}_{t_1}^{(i_1)}\ldots
d{\bf w}_{t_k}^{(i_k)},
\end{equation}

\vspace{2mm}
\noindent
where every $\psi_l(\tau)$ $(l=1,\ldots,k)$ is a non-random
function 
on $[t,T],$ ${\bf w}_{\tau}^{(i)}={\bf f}_{\tau}^{(i)}$
for $i=1,\ldots,m$ and
${\bf w}_{\tau}^{(0)}=\tau,$\ \
$i_1,\ldots,i_k = 0, 1,\ldots,m,$

\vspace{-2mm}
$$
\int\limits\ \hbox{and}\ \int\limits^{*}
$$ 

\vspace{2mm}
\noindent
denote Ito and 
Stratonovich stochastic integrals,
respectively (in this paper, 
we use the definition of the Stratonovich stochastic integral from \cite{KlPl2}).

Note that $\psi_l(\tau)\equiv 1$ $(l=1,\ldots,k)$ and
$i_1,\ldots,i_k = 0, 1,\ldots,m$ in  
\cite{KlPl2}-\cite{KlPl1}. At the same time
$\psi_l(\tau)\equiv (t-\tau)^{q_l}$ ($l=1,\ldots,k$; 
$q_1,\ldots,q_k=0, 1, 2,\ldots $) and  $i_1,\ldots,i_k = 1,\ldots,m$ in
\cite{kk5}-\cite{arxiv-24999}.
 
We want to mention in short that there are 
two main criteria of numerical methods convergence for Ito stochastic
differential equations:  a strong or 
mean-square
criterion and a 
weak criterion, where the subject of approximation is not the solution 
of Ito stochastic differential equation, simply stated, but the 
distribution of Ito stochastic differential equation solution \cite{KlPl2}.
Both of the above criteria are independent, that is,
generally speaking,  
the fulfillment of a strong criterion does not imply
the fulfillment of a weak criterion, and vice versa.
Each of two convergence criteria is oriented on solution of specific 
classes of mathematical problems connected with stochastic differential 
equations.
   
Using the strong numerical methods, we may build 
sample pathes
of Ito stochastic differential equation numerically. These 
methods require the combined mean-square approximation of collections 
of iterated Ito and Stratonovich stochastic integrals. Effective solution 
of this problem composes one of the subjects of this article.

The strong numerical methods are used for constructing new mathematical 
models on the basis of Ito stochastic differential equations and 
also for solving some mathematical problems connected with Ito stochastic 
differential 
equations. Among this problems we mention the following:
signal filtering under the influence
of random disturbances in various statements, 
stochastic 
optimal control,  
testing estimation procedures of parameters of stochastic 
systems \cite{KlPl2}.

\vspace{5mm}

\section{Explicit One-Step Strong Numerical Schemes of Orders 2.0 and 
2.5 Based on the Unified Taylor--Ito expansion}

\vspace{5mm}

Consider the partition $\{\tau_j\}_{j=0}^N$ 
of the interval $[0, T]$ such that

\vspace{1mm}
$$
0=\tau_0<\ldots <\tau_N=T,\ \ \
\Delta_N=
\hbox{\vtop{\offinterlineskip\halign{
\hfil#\hfil\cr
{\rm max}\cr
$\stackrel{}{{}_{0\le j\le N-1}}$\cr
}} }\Delta\tau_j,\ \ \ \Delta\tau_j=\tau_{j+1}-\tau_j.
$$

\vspace{4mm}

Let ${\bf y}_{\tau_j}\stackrel{\sf def}{=}
{\bf y}_{j},$\ $j=0, 1,\ldots,N$ be a time discrete approximation
of the process ${\bf x}_t,$ $t\in[0,T],$ which is a solution of the Ito
stochastic differential equation (\ref{1.5.2}). 

\vspace{2mm}

{\bf Definiton 1}\ \cite{KlPl2}.\
{\it We will say that a time discrete approximation 
${\bf y}_{j},$\ $j=0, 1,\ldots,N,$
corresponding to the maximal step of discretization $\Delta_N,$
converges strongly with order
$\gamma>0$ at time moment 
$T$ to the process ${\bf x}_t,$ $t\in[0,T]$,
if there exists a constant $C>0,$ which does not depend on 
$\Delta_N,$ and a $\delta>0$ such that 

\vspace{-1mm}
$$
{\sf M}\{|{\bf x}_T-{\bf y}_T|\}\le
C(\Delta_N)^{\gamma}
$$

\vspace{2mm}
\noindent
for each $\Delta_N\in(0, \delta).$}

\vspace{2mm}

Consider the explicit one-step strong numerical scheme of order 2.5
based on the so-called unified Taylor--Ito expansion 
\cite{2006}, \cite{2017-1}-\cite{2010-1}

\vspace{2mm}
$$
{\bf y}_{p+1}={\bf y}_p+\sum_{i_{1}=1}^{m}B_{i_{1}}
\hat I_{(0)\tau_{p+1},\tau_p}^{(i_{1})}+\Delta{\bf a}
+\sum_{i_{1},i_{2}=1}^{m}G_{i_{2}}
B_{i_{1}}\hat I_{(00)\tau_{p+1},\tau_p}^{(i_{2}i_{1})}+
$$

\vspace{1mm}
$$
+
\sum_{i_{1}=1}^{m}\Biggl(G_{i_{1}}{\bf a}\left(
\Delta \hat I_{(0)\tau_{p+1},\tau_p}^{(i_{1})}+
\hat I_{(1)\tau_{p+1},\tau_p}^{(i_{1})}\right)
-LB_{i_{1}}\hat I_{(1)\tau_{p+1},\tau_p}^{(i_{1})}\Biggr)+
$$

\vspace{1mm}
$$
+\sum_{i_{1},i_{2},i_{3}=1}^{m} G_{i_{3}}G_{i_{2}}
B_{i_{1}} \hat I_{(000)\tau_{p+1},\tau_p}^{(i_{3}i_{2}i_{1})}+
\frac{\Delta^2}{2}L{\bf a}+
$$

\vspace{1mm}
$$
+\sum_{i_{1},i_{2}=1}^{m}
\Biggl(G_{i_{2}}LB_{i_{1}}\left(
\hat I_{(10)\tau_{p+1},\tau_p}^{(i_{2}i_{1})}-
\hat I_{(01)\tau_{p+1},\tau_p}^{(i_{2}i_{1})}
\right)
-LG_{i_{2}}B_{i_{1}}\hat I_{(10)\tau_{p+1},\tau_p}^{(i_{2}i_{1})}
+\Biggr.
$$

\vspace{1mm}
$$
\Biggl.+G_{i_{2}}G_{i_{1}}{\bf a}\left(
\hat I_{(01)\tau_{p+1},\tau_p}^{(i_{2}i_{1})}+
\Delta \hat I_{(00)\tau_{p+1},\tau_p}^{(i_{2}i_{1})}
\right)\Biggr)+
$$

\vspace{1mm}
$$
+
\sum_{i_{1},i_{2},i_{3},i_{4}=1}^{m}G_{i_{4}}G_{i_{3}}G_{i_{2}}
B_{i_{1}}\hat I_{(0000)\tau_{p+1},\tau_p}^{(i_{4}i_{3}i_{2}i_{1})}+
\frac{\Delta^3}{6}LL{\bf a}+
$$

\vspace{1mm}
$$
+\sum_{i_{1}=1}^{m}\Biggl(G_{i_{1}}L{\bf a}\left(\frac{1}{2}
\hat I_{(2)\tau_{p+1},\tau_p}
^{(i_{1})}+\Delta \hat I_{(1)\tau_{p+1},\tau_p}^{(i_{1})}+
\frac{\Delta^2}{2}\hat I_{(0)\tau_{p+1},\tau_p}^{(i_{1})}\right)\Biggr.+
$$

\vspace{1mm}
$$
\Biggl.+\frac{1}{2}LL B_{i_{1}}\hat I_{(2)\tau_{p+1},\tau_p}^{(i_{1})}-
LG_{i_{1}}{\bf a}\left(\hat I_{(2)\tau_{p+1},\tau_p}^{(i_{1})}+
\Delta \hat I_{(1)\tau_{p+1},\tau_p}^{(i_{1})}\right)\Biggr)+
$$

\vspace{1mm}
$$
+
\sum_{i_{1},i_{2},i_{3}=1}^m\Biggl(
G_{i_{3}}LG_{i_{2}}B_{i_{1}}
\left(\hat I_{(100)\tau_{p+1},\tau_p}
^{(i_{3}i_{2}i_{1})}-\hat I_{(010)\tau_{p+1},\tau_p}
^{(i_{3}i_{2}i_{1})}\right)
\Biggr.+
$$

\vspace{1mm}
$$
+G_{i_{3}}G_{i_{2}}LB_{i_{1}}\left(
\hat I_{(010)\tau_{p+1},\tau_p}^{(i_{3}i_{2}i_{1})}-
\hat I_{(001)\tau_{p+1},\tau_p}^{(i_{3}i_{2}i_{1})}\right)+
$$

\vspace{1mm}

$$
+
G_{i_{3}}G_{i_{2}}G_{i_{1}}{\bf a}
\left(\Delta \hat I_{(000)\tau_{p+1},\tau_p}^{(i_{3}i_{2}i_{1})}+
\hat I_{(001)\tau_{p+1},\tau_p}^{(i_{3}i_{2}i_{1})}\right)
-
$$

\vspace{1mm}
$$
\Biggl.-LG_{i_{3}}G_{i_{2}}B_{i_{1}}
\hat I_{(100)\tau_{p+1},\tau_p}^{(i_{3}i_{2}i_{1})}\Biggr)+
$$

\vspace{1mm}
\begin{equation}
\label{4.45}
+\sum_{i_{1},i_{2},i_{3},i_{4},i_{5}=1}^m
G_{i_{5}}G_{i_{4}}G_{i_{3}}G_{i_{2}}B_{i_{1}}
\hat I_{(00000)\tau_{p+1},\tau_p}^{(i_{5}i_{4}i_{3}i_{2}i_{1})},
\end{equation}

\vspace{5mm}
\noindent
where $\Delta=T/N$ $(N>1)$ is a constant (for simplicity) 
integration step,\
$\tau_p=p\Delta$ $(p=0, 1,\ldots,N)$,\
$\hat I_{(l_1\ldots\hspace{0.2mm} l_k)s,t}^{(i_1\ldots i_k)}$ 
denotes
an approximation of the iterated Ito stochastic integral of multiplicity
$k$

\vspace{-1mm}
\begin{equation}
\label{ll1}
I_{(l_1\ldots\hspace{0.2mm} l_k)s,t}^{(i_1\ldots i_k)}=
 \int\limits^ {s} _ {t} (t-\tau _
{k}) ^ {l_ {k}} 
\ldots \int\limits^ {\tau _ {2}} _ {t} (t-\tau _ {1}) ^ {l_ {1}} d
{\bf f} ^ {(i_ {1})} _ {\tau_ {1}} \ldots 
d {\bf f} _ {\tau_ {k}} ^ {(i_ {k})},
\end{equation}

\vspace{3mm}
$$
L= {\partial \over \partial t}
+ \sum^ {n} _ {i=1} {\bf a}_i ({\bf x},  t) 
{\partial  \over  \partial  {\bf  x}_i}
+ {1\over 2} \sum^ {m} _ {j=1} \sum^ {n} _ {l,i=1}
B_{lj} ({\bf x}, t) B_{ij} ({\bf x}, t) {\partial
^{2} \over \partial {\bf x}_l \partial {\bf x}_i},
$$

\vspace{3mm}
$$
G_i = \sum^ {n} _ {j=1} B_{ji} ({\bf x}, t)
{\partial  \over \partial {\bf x}_j},\ \ \ 
i=1,\ldots,m,
$$

\vspace{3mm}
\noindent
$l_1,\ldots, l_k=0, 1, 2,$\ \
$i_1,\ldots, i_k=1,\ldots,m,$\ \ $k=1, 2,\ldots, 5$,\ \
$B_i$ and $B_{ij}$ are the $i$th column and the $ij$th
component of the matrix function $B$,
${\bf a}_i$ is the $i$th component of the vector function ${\bf a},$
${\bf x}_i$ is the $i$th component
of the column ${\bf x}$, 
the functions  

\vspace{-1mm}
$$
B_{i_{1}},\ {\bf a},\ G_{i_{2}}B_{i_{1}},\
G_{i_{1}}{\bf a},\ LB_{i_{1}},\ G_{i_{3}}G_{i_{2}}B_{i_{1}},\ 
L{\bf a},\ LL{\bf a},\
G_{i_{2}}LB_{i_{1}},
$$

\vspace{-3mm}
$$
LG_{i_{2}}B_{i_{1}},\ G_{i_{2}}G_{i_{1}}{\bf a},\
G_{i_{4}}G_{i_{3}}G_{i_{2}}B_{i_{1}},\ G_{i_{1}}L{\bf a},\
LLB_{i_{1}},\ LG_{i_{1}}{\bf a},\ G_{i_{3}}LG_{i_{2}}B_{i_{1}},\
G_{i_{3}}G_{i_{2}}LB_{i_{1}},
$$

\vspace{-3mm}
$$
G_{i_{3}}G_{i_{2}}G_{i_{1}}{\bf a},\
LG_{i_{3}}G_{i_{2}}B_{i_{1}},\ 
G_{i_{5}}G_{i_{4}}G_{i_{3}}G_{i_{2}}B_{i_{1}}
$$

\vspace{3mm}
\noindent
are calculated at the point $({\bf y}_p,p).$

Under the standard conditions \cite{KlPl2}, \cite{2006} the numerical 
scheme (\ref{4.45}) has 
strong order 2.5 of convergence. 
The major emphasis below will be placed on the 
approximation of the iterated
Ito stochastic integrals appearing in (\ref{4.45}). Therefore, among 
the mentioned conditions, we note only the
approximation condi\-ti\-on for iterated Ito
stochastic integrals \cite{KlPl2}, \cite{2006}, 
which has the form

\vspace{-1mm}
\begin{equation}
\label{4.3}
{\sf M}\left\{\biggl(I_{(l_{1}\ldots\hspace{0.2mm} l_{k})\tau_{p+1},\tau_p}
^{(i_{1}\ldots i_{k})} 
- \hat I_{(l_{1}\ldots\hspace{0.2mm} l_{k})\tau_{p+1},
\tau_p}^{(i_{1}\ldots i_{k})}
\biggr)^2\right\}\le C\Delta^{6},
\end{equation}

\vspace{2mm}
\noindent
where constant $C$ is independent of
$\Delta$.

Note that if we exclude from (\ref{4.45}) the terms starting from the
term $\Delta^3 LL{\bf a}/6$, then we will have the explicit 
one-step strong numerical scheme of order 2.0 \cite{KlPl2}, 
\cite{2006}, \cite{2017-1}-\cite{2010-1}.

Using the numerical scheme (\ref{4.45}) or its modifications based 
on the Taylor--Ito expansion \cite{KlPl1},
the implicit or multistep analogues of (\ref{4.45}) can be constructed
\cite{KlPl2}, \cite{2006}, \cite{2017-1}-\cite{2010-1}. The set of the
iterated Ito stochastic integrals to be approximated for implementing 
these modifications is the same
as for the numerical scheme (\ref{4.45}) itself.
Interestingly, the truncated unified Taylor--Ito expansion (the 
foundation of the numerical
scheme (\ref{4.45})) contains 12 different types of iterated Ito 
stochastic integrals 
of the form (\ref{ll1}), which cannot be
interconnected by linear relations \cite{2006}, \cite{2017-1}-\cite{2010-1}. 
The analogous 
Taylor--Ito expansion \cite{KlPl2}, \cite{KlPl1} contains
17 different types of iterated Ito stochastic integrals, part of which 
are interconnected by linear relations
and part of which have a higher multiplicity than the iterated 
Ito stochastic integrals (\ref{ll1}). This
fact well explains the use of the numerical scheme (\ref{4.45}).

One of the main problems arising in the implementation of the 
numerical scheme (\ref{4.45}) is the joint
numerical modeling of the iterated Ito stochastic integrals 
figuring in (\ref{4.45}). In the subsequent sections,
we will consider an efficient numerical modeling method for 
the iterated Ito stochastic integrals
and also demonstrate which stochastic integrals (Ito or Stratonovich) 
are preferable for numerical
modeling with a correct estimation of the mean-square approximation error.

\vspace{5mm}

\section{Method of Numerical Modeling for Iterated Ito 
Stochastic Integrals
Based on Genegalized Multiple Fourier Series. Direct Approach}

\vspace{5mm}

An efficient numerical modeling method for the iterated Ito 
stochastic integrals (\ref{ito}) based on generalized multiple
Fourier series was proposed in \cite{2006} (2006); also see 
\cite{2011-2}-\cite{new-art-1xxys}. 
This method rests on an important
result presented below.

Suppose that every $\psi_l(\tau)$ $(l=1,\ldots,k)$ is a 
non-random function from the space $L_2([t, T])$.
Define the following function on the hypercube $[t, T]^k$

\vspace{-1mm}
\begin{equation}
\label{ppp}
K(t_1,\ldots,t_k)=
\begin{cases}
\psi_1(t_1)\ldots \psi_k(t_k)\ &\hbox{for}\ \ t_1<\ldots<t_k\\
~\\
~\\
0\ &\hbox{otherwise}
\end{cases},\ \ \ \ t_1,\ldots,t_k\in[t, T],\ \ \ \ k\ge 2,
\end{equation}

\vspace{3mm}
\noindent
and 
$K(t_1)\equiv\psi_1(t_1)$ for $t_1\in[t, T].$

Suppose that $\{\phi_j(x)\}_{j=0}^{\infty}$
is a complete orthonormal system of functions in the space
$L_2([t, T])$. 
The function $K(t_1,\ldots,t_k)$ belongs to the 
space $L_2([t, T]^k).$
At this situation it is well known that the generalized 
multiple Fourier series 
of $K(t_1,\ldots,t_k)\in L_2([t, T]^k)$ is converging 
to $K(t_1,\ldots,t_k)$ in the hypercube $[t, T]^k$ in 
the mean-square sense, i.e.

\vspace{-1mm}
$$
\hbox{\vtop{\offinterlineskip\halign{
\hfil#\hfil\cr
{\rm lim}\cr
$\stackrel{}{{}_{p_1,\ldots,p_k\to \infty}}$\cr
}} }\left\Vert
K(t_1,\ldots,t_k)-
\sum_{j_1=0}^{p_1}\ldots \sum_{j_k=0}^{p_k}
C_{j_k\ldots j_1}\prod_{l=1}^{k} \phi_{j_l}(t_l)
\right\Vert_{L_2([t, T]^k)}=0,
$$

\vspace{3mm}
\noindent
where
\begin{equation}
\label{ppppa}
C_{j_k\ldots j_1}=\int\limits_{[t,T]^k}
K(t_1,\ldots,t_k)\prod_{l=1}^{k}\phi_{j_l}(t_l)dt_1\ldots dt_k,
\end{equation}

\vspace{2mm}
$$
\left\Vert f\right\Vert_{L_2([t, T]^k)}=\left(\int\limits_{[t,T]^k}
f^2(t_1,\ldots,t_k)dt_1\ldots dt_k\right)^{1/2}.
$$

\vspace{5mm}

Consider the partition $\{\tau_j\}_{j=0}^N$ of $[t,T]$ such that

\vspace{1mm}
\begin{equation}
\label{1111}
t=\tau_0<\ldots <\tau_N=T,\ \ \
\Delta_N=
\hbox{\vtop{\offinterlineskip\halign{
\hfil#\hfil\cr
{\rm max}\cr
$\stackrel{}{{}_{0\le j\le N-1}}$\cr
}} }\Delta\tau_j\to 0\ \ \hbox{if}\ \ N\to \infty,\ \ \
\Delta\tau_j=\tau_{j+1}-\tau_j.
\end{equation}

\vspace{5mm}

{\bf Theorem 1}\ \cite{2006}-\cite{200a}, \cite{301a}-\cite{arxiv-12}, 
\cite{arxiv-24}-\cite{new-art-1xxys}.
{\it Suppose that
every $\psi_l(\tau)$ $(l=1,\ldots, k)$ is a 
continuous non-random func\-tion on 
$[t, T]$ and
$\{\phi_j(x)\}_{j=0}^{\infty}$ is a complete orthonormal system  
of continuous func\-ti\-ons in the space $L_2([t,T]).$ Then

$$
J[\psi^{(k)}]_{T,t}\  =\ 
\hbox{\vtop{\offinterlineskip\halign{
\hfil#\hfil\cr
{\rm l.i.m.}\cr
$\stackrel{}{{}_{p_1,\ldots,p_k\to \infty}}$\cr
}} }\sum_{j_1=0}^{p_1}\ldots\sum_{j_k=0}^{p_k}
C_{j_k\ldots j_1}\Biggl(
\prod_{l=1}^k\zeta_{j_l}^{(i_l)}\ -
\Biggr.
$$

\vspace{1mm}
\begin{equation}
\label{tyyy}
-\ \Biggl.
\hbox{\vtop{\offinterlineskip\halign{
\hfil#\hfil\cr
{\rm l.i.m.}\cr
$\stackrel{}{{}_{N\to \infty}}$\cr
}} }\sum_{(l_1,\ldots,l_k)\in {\rm G}_k}
\phi_{j_{1}}(\tau_{l_1})
\Delta{\bf w}_{\tau_{l_1}}^{(i_1)}\ldots
\phi_{j_{k}}(\tau_{l_k})
\Delta{\bf w}_{\tau_{l_k}}^{(i_k)}\Biggr),
\end{equation}

\vspace{5mm}
\noindent
where $J[\psi^{(k)}]_{T,t}$ is defined by {\rm (\ref{ito}),}

$$
{\rm G}_k={\rm H}_k\backslash{\rm L}_k,\ \ \
{\rm H}_k=\{(l_1,\ldots,l_k):\ l_1,\ldots,l_k=0,\ 1,\ldots,N-1\},
$$

\vspace{1mm}
$$
{\rm L}_k=\{(l_1,\ldots,l_k):\ l_1,\ldots,l_k=0,\ 1,\ldots,N-1;\
l_g\ne l_r\ (g\ne r);\ g, r=1,\ldots,k\},
$$

\vspace{4mm}
\noindent
${\rm l.i.m.}$ is a limit in the mean-square sense$,$
$i_1,\ldots,i_k=0,1,\ldots,m,$

\vspace{-1mm}
\begin{equation}
\label{rr23}
\zeta_{j}^{(i)}=
\int\limits_t^T \phi_{j}(s) d{\bf w}_s^{(i)}
\end{equation} 

\vspace{2mm}
\noindent
are independent standard Gaussian random variables
for various
$i$ or $j$ {\rm(}if $i\ne 0${\rm),}
$C_{j_k\ldots j_1}$ is the Fourier coefficient {\rm(\ref{ppppa}),}
$\Delta{\bf w}_{\tau_{j}}^{(i)}=
{\bf w}_{\tau_{j+1}}^{(i)}-{\bf w}_{\tau_{j}}^{(i)}$
$(i=0, 1,\ldots,m),$
$\left\{\tau_{j}\right\}_{j=0}^{N}$ is a partition of
the interval $[t, T],$ which satisfies the condition {\rm (\ref{1111})}.
}

\vspace{2mm}

The convergence in the mean of degree 
$2n$ $(n\in \mathbb{N})$ \cite{2018a}-\cite{2013}
as well as the convergence
with probability 1 \cite{2018a}-\cite{2018aaa}, \cite{arxiv-1},
\cite{arxiv-3} are proved for the approximations from Theorem 1.

Moreover, the complete orthonormal systems of Haar and 
Rademacher--Walsh functions in the space $L_2([t,T])$ 
can also be applied in Theorem 1
\cite{2006}-\cite{2013}.
The modification of Theorem 1 for 
complete orthonormal with weigth $r(x)\ge 0$ systems
of functions in the space $L_2([t,T])$ can be found in 
\cite{2018}-\cite{2018aaa}, \cite{arxiv-26b}.

In order to evaluate the significance of Theorem 1 for practice we will
demonstrate its transfor\-med particular cases for 
$k=1,\ldots,5$ 
\cite{2006}-\cite{200a}, \cite{301a}-\cite{arxiv-12}, 
\cite{arxiv-24}-\cite{new-art-1xxys}

\begin{equation}
\label{a1}
J[\psi^{(1)}]_{T,t}
=\hbox{\vtop{\offinterlineskip\halign{
\hfil#\hfil\cr
{\rm l.i.m.}\cr
$\stackrel{}{{}_{p_1\to \infty}}$\cr
}} }\sum_{j_1=0}^{p_1}
C_{j_1}\zeta_{j_1}^{(i_1)},
\end{equation}

\vspace{2mm}
\begin{equation}
\label{a2}
J[\psi^{(2)}]_{T,t}
=\hbox{\vtop{\offinterlineskip\halign{
\hfil#\hfil\cr
{\rm l.i.m.}\cr
$\stackrel{}{{}_{p_1,p_2\to \infty}}$\cr
}} }\sum_{j_1=0}^{p_1}\sum_{j_2=0}^{p_2}
C_{j_2j_1}\Biggl(\zeta_{j_1}^{(i_1)}\zeta_{j_2}^{(i_2)}
-{\bf 1}_{\{i_1=i_2\ne 0\}}
{\bf 1}_{\{j_1=j_2\}}\Biggr),
\end{equation}

\vspace{6mm}

$$
J[\psi^{(3)}]_{T,t}=
\hbox{\vtop{\offinterlineskip\halign{
\hfil#\hfil\cr
{\rm l.i.m.}\cr
$\stackrel{}{{}_{p_1,p_2,p_3\to \infty}}$\cr
}} }\sum_{j_1=0}^{p_1}\sum_{j_2=0}^{p_2}\sum_{j_3=0}^{p_3}
C_{j_3j_2j_1}\Biggl(
\zeta_{j_1}^{(i_1)}\zeta_{j_2}^{(i_2)}\zeta_{j_3}^{(i_3)}
-\Biggr.
$$
\begin{equation}
\label{a3}
-\Biggl.
{\bf 1}_{\{i_1=i_2\ne 0\}}
{\bf 1}_{\{j_1=j_2\}}
\zeta_{j_3}^{(i_3)}
-{\bf 1}_{\{i_2=i_3\ne 0\}}
{\bf 1}_{\{j_2=j_3\}}
\zeta_{j_1}^{(i_1)}-
{\bf 1}_{\{i_1=i_3\ne 0\}}
{\bf 1}_{\{j_1=j_3\}}
\zeta_{j_2}^{(i_2)}\Biggr),
\end{equation}

\vspace{6mm}

$$
J[\psi^{(4)}]_{T,t}
=
\hbox{\vtop{\offinterlineskip\halign{
\hfil#\hfil\cr
{\rm l.i.m.}\cr
$\stackrel{}{{}_{p_1,\ldots,p_4\to \infty}}$\cr
}} }\sum_{j_1=0}^{p_1}\ldots\sum_{j_4=0}^{p_4}
C_{j_4\ldots j_1}\Biggl(
\prod_{l=1}^4\zeta_{j_l}^{(i_l)}
\Biggr.
-
$$
$$
-
{\bf 1}_{\{i_1=i_2\ne 0\}}
{\bf 1}_{\{j_1=j_2\}}
\zeta_{j_3}^{(i_3)}
\zeta_{j_4}^{(i_4)}
-
{\bf 1}_{\{i_1=i_3\ne 0\}}
{\bf 1}_{\{j_1=j_3\}}
\zeta_{j_2}^{(i_2)}
\zeta_{j_4}^{(i_4)}-
$$
$$
-
{\bf 1}_{\{i_1=i_4\ne 0\}}
{\bf 1}_{\{j_1=j_4\}}
\zeta_{j_2}^{(i_2)}
\zeta_{j_3}^{(i_3)}
-
{\bf 1}_{\{i_2=i_3\ne 0\}}
{\bf 1}_{\{j_2=j_3\}}
\zeta_{j_1}^{(i_1)}
\zeta_{j_4}^{(i_4)}-
$$
$$
-
{\bf 1}_{\{i_2=i_4\ne 0\}}
{\bf 1}_{\{j_2=j_4\}}
\zeta_{j_1}^{(i_1)}
\zeta_{j_3}^{(i_3)}
-
{\bf 1}_{\{i_3=i_4\ne 0\}}
{\bf 1}_{\{j_3=j_4\}}
\zeta_{j_1}^{(i_1)}
\zeta_{j_2}^{(i_2)}+
$$
$$
+
{\bf 1}_{\{i_1=i_2\ne 0\}}
{\bf 1}_{\{j_1=j_2\}}
{\bf 1}_{\{i_3=i_4\ne 0\}}
{\bf 1}_{\{j_3=j_4\}}
+
{\bf 1}_{\{i_1=i_3\ne 0\}}
{\bf 1}_{\{j_1=j_3\}}
{\bf 1}_{\{i_2=i_4\ne 0\}}
{\bf 1}_{\{j_2=j_4\}}+
$$
\begin{equation}
\label{a4}
+\Biggl.
{\bf 1}_{\{i_1=i_4\ne 0\}}
{\bf 1}_{\{j_1=j_4\}}
{\bf 1}_{\{i_2=i_3\ne 0\}}
{\bf 1}_{\{j_2=j_3\}}\Biggr),
\end{equation}

\vspace{6mm}

$$
J[\psi^{(5)}]_{T,t}
=\hbox{\vtop{\offinterlineskip\halign{
\hfil#\hfil\cr
{\rm l.i.m.}\cr
$\stackrel{}{{}_{p_1,\ldots,p_5\to \infty}}$\cr
}} }\sum_{j_1=0}^{p_1}\ldots\sum_{j_5=0}^{p_5}
C_{j_5\ldots j_1}\Biggl(
\prod_{l=1}^5\zeta_{j_l}^{(i_l)}
-\Biggr.
$$
$$
-
{\bf 1}_{\{i_1=i_2\ne 0\}}
{\bf 1}_{\{j_1=j_2\}}
\zeta_{j_3}^{(i_3)}
\zeta_{j_4}^{(i_4)}
\zeta_{j_5}^{(i_5)}-
{\bf 1}_{\{i_1=i_3\ne 0\}}
{\bf 1}_{\{j_1=j_3\}}
\zeta_{j_2}^{(i_2)}
\zeta_{j_4}^{(i_4)}
\zeta_{j_5}^{(i_5)}-
$$
$$
-
{\bf 1}_{\{i_1=i_4\ne 0\}}
{\bf 1}_{\{j_1=j_4\}}
\zeta_{j_2}^{(i_2)}
\zeta_{j_3}^{(i_3)}
\zeta_{j_5}^{(i_5)}-
{\bf 1}_{\{i_1=i_5\ne 0\}}
{\bf 1}_{\{j_1=j_5\}}
\zeta_{j_2}^{(i_2)}
\zeta_{j_3}^{(i_3)}
\zeta_{j_4}^{(i_4)}-
$$
$$
-
{\bf 1}_{\{i_2=i_3\ne 0\}}
{\bf 1}_{\{j_2=j_3\}}
\zeta_{j_1}^{(i_1)}
\zeta_{j_4}^{(i_4)}
\zeta_{j_5}^{(i_5)}-
{\bf 1}_{\{i_2=i_4\ne 0\}}
{\bf 1}_{\{j_2=j_4\}}
\zeta_{j_1}^{(i_1)}
\zeta_{j_3}^{(i_3)}
\zeta_{j_5}^{(i_5)}-
$$
$$
-
{\bf 1}_{\{i_2=i_5\ne 0\}}
{\bf 1}_{\{j_2=j_5\}}
\zeta_{j_1}^{(i_1)}
\zeta_{j_3}^{(i_3)}
\zeta_{j_4}^{(i_4)}
-{\bf 1}_{\{i_3=i_4\ne 0\}}
{\bf 1}_{\{j_3=j_4\}}
\zeta_{j_1}^{(i_1)}
\zeta_{j_2}^{(i_2)}
\zeta_{j_5}^{(i_5)}-
$$
$$
-
{\bf 1}_{\{i_3=i_5\ne 0\}}
{\bf 1}_{\{j_3=j_5\}}
\zeta_{j_1}^{(i_1)}
\zeta_{j_2}^{(i_2)}
\zeta_{j_4}^{(i_4)}
-{\bf 1}_{\{i_4=i_5\ne 0\}}
{\bf 1}_{\{j_4=j_5\}}
\zeta_{j_1}^{(i_1)}
\zeta_{j_2}^{(i_2)}
\zeta_{j_3}^{(i_3)}+
$$
$$
+
{\bf 1}_{\{i_1=i_2\ne 0\}}
{\bf 1}_{\{j_1=j_2\}}
{\bf 1}_{\{i_3=i_4\ne 0\}}
{\bf 1}_{\{j_3=j_4\}}\zeta_{j_5}^{(i_5)}+
{\bf 1}_{\{i_1=i_2\ne 0\}}
{\bf 1}_{\{j_1=j_2\}}
{\bf 1}_{\{i_3=i_5\ne 0\}}
{\bf 1}_{\{j_3=j_5\}}\zeta_{j_4}^{(i_4)}+
$$
$$
+
{\bf 1}_{\{i_1=i_2\ne 0\}}
{\bf 1}_{\{j_1=j_2\}}
{\bf 1}_{\{i_4=i_5\ne 0\}}
{\bf 1}_{\{j_4=j_5\}}\zeta_{j_3}^{(i_3)}+
{\bf 1}_{\{i_1=i_3\ne 0\}}
{\bf 1}_{\{j_1=j_3\}}
{\bf 1}_{\{i_2=i_4\ne 0\}}
{\bf 1}_{\{j_2=j_4\}}\zeta_{j_5}^{(i_5)}+
$$
$$
+
{\bf 1}_{\{i_1=i_3\ne 0\}}
{\bf 1}_{\{j_1=j_3\}}
{\bf 1}_{\{i_2=i_5\ne 0\}}
{\bf 1}_{\{j_2=j_5\}}\zeta_{j_4}^{(i_4)}+
{\bf 1}_{\{i_1=i_3\ne 0\}}
{\bf 1}_{\{j_1=j_3\}}
{\bf 1}_{\{i_4=i_5\ne 0\}}
{\bf 1}_{\{j_4=j_5\}}\zeta_{j_2}^{(i_2)}+
$$
$$
+
{\bf 1}_{\{i_1=i_4\ne 0\}}
{\bf 1}_{\{j_1=j_4\}}
{\bf 1}_{\{i_2=i_3\ne 0\}}
{\bf 1}_{\{j_2=j_3\}}\zeta_{j_5}^{(i_5)}+
{\bf 1}_{\{i_1=i_4\ne 0\}}
{\bf 1}_{\{j_1=j_4\}}
{\bf 1}_{\{i_2=i_5\ne 0\}}
{\bf 1}_{\{j_2=j_5\}}\zeta_{j_3}^{(i_3)}+
$$
$$
+
{\bf 1}_{\{i_1=i_4\ne 0\}}
{\bf 1}_{\{j_1=j_4\}}
{\bf 1}_{\{i_3=i_5\ne 0\}}
{\bf 1}_{\{j_3=j_5\}}\zeta_{j_2}^{(i_2)}+
{\bf 1}_{\{i_1=i_5\ne 0\}}
{\bf 1}_{\{j_1=j_5\}}
{\bf 1}_{\{i_2=i_3\ne 0\}}
{\bf 1}_{\{j_2=j_3\}}\zeta_{j_4}^{(i_4)}+
$$
$$
+
{\bf 1}_{\{i_1=i_5\ne 0\}}
{\bf 1}_{\{j_1=j_5\}}
{\bf 1}_{\{i_2=i_4\ne 0\}}
{\bf 1}_{\{j_2=j_4\}}\zeta_{j_3}^{(i_3)}+
{\bf 1}_{\{i_1=i_5\ne 0\}}
{\bf 1}_{\{j_1=j_5\}}
{\bf 1}_{\{i_3=i_4\ne 0\}}
{\bf 1}_{\{j_3=j_4\}}\zeta_{j_2}^{(i_2)}+
$$
$$
+
{\bf 1}_{\{i_2=i_3\ne 0\}}
{\bf 1}_{\{j_2=j_3\}}
{\bf 1}_{\{i_4=i_5\ne 0\}}
{\bf 1}_{\{j_4=j_5\}}\zeta_{j_1}^{(i_1)}+
{\bf 1}_{\{i_2=i_4\ne 0\}}
{\bf 1}_{\{j_2=j_4\}}
{\bf 1}_{\{i_3=i_5\ne 0\}}
{\bf 1}_{\{j_3=j_5\}}\zeta_{j_1}^{(i_1)}+
$$
\begin{equation}
\label{a5}
+\Biggl.
{\bf 1}_{\{i_2=i_5\ne 0\}}
{\bf 1}_{\{j_2=j_5\}}
{\bf 1}_{\{i_3=i_4\ne 0\}}
{\bf 1}_{\{j_3=j_4\}}\zeta_{j_1}^{(i_1)}\Biggr),
\end{equation}

\vspace{6mm}
\noindent
where ${\bf 1}_A$ is the indicator of the set $A$.

We will consider the case $i_1,\ldots,i_5=1,\ldots,m$.
Obviously, this case corresponds to the numerical method (\ref{4.45}).

For further consideration, let us 
consider the generalization of formulas (\ref{a1})--(\ref{a5})                 
for the case of an arbitrary multiplicity $k$ $(k\in\mathbb{N})$ of 
the iterated Ito stochastic integral $J[\psi^{(k)}]_{T,t}$ defined by (\ref{ito}).
In order to do this, let us
introduce some notations. 
Consider the unordered
set $\{1, 2, \ldots, k\}$ 
and separate it into two parts:
the first part consists of $r$ unordered 
pairs (sequence order of these pairs is also unimportant) and the 
second one consists of the 
remaining $k-2r$ numbers.
So, we have

\begin{equation}
\label{leto5007}
(\{
\underbrace{\{g_1, g_2\}, \ldots, 
\{g_{2r-1}, g_{2r}\}}_{\small{\hbox{part 1}}}
\},
\{\underbrace{q_1, \ldots, q_{k-2r}}_{\small{\hbox{part 2}}}
\}),
\end{equation}

\vspace{4mm}
\noindent
where 

\vspace{-2mm}
$$
\{g_1, g_2, \ldots, 
g_{2r-1}, g_{2r}, q_1, \ldots, q_{k-2r}\}=\{1, 2, \ldots, k\},
$$

\vspace{4mm}
\noindent
braces   
mean an unordered 
set, and pa\-ren\-the\-ses mean an ordered set.

We will say that (\ref{leto5007}) is a partition 
and consider the sum with respect to all possible
partitions

\begin{equation}
\label{leto5008}
\sum_{\stackrel{(\{\{g_1, g_2\}, \ldots, 
\{g_{2r-1}, g_{2r}\}\}, \{q_1, \ldots, q_{k-2r}\})}
{{}_{\{g_1, g_2, \ldots, 
g_{2r-1}, g_{2r}, q_1, \ldots, q_{k-2r}\}=\{1, 2, \ldots, k\}}}}
a_{g_1 g_2, \ldots, 
g_{2r-1} g_{2r}, q_1 \ldots q_{k-2r}}.
\end{equation}

\vspace{4mm}

Below there are several examples of sums in the form (\ref{leto5008})

\vspace{2mm}
$$
\sum_{\stackrel{(\{g_1, g_2\})}{{}_{\{g_1, g_2\}=\{1, 2\}}}}
a_{g_1 g_2}=a_{12},
$$

\vspace{3mm}
$$
\sum_{\stackrel{(\{\{g_1, g_2\}, \{g_3, g_4\}\})}
{{}_{\{g_1, g_2, g_3, g_4\}=\{1, 2, 3, 4\}}}}
a_{g_1 g_2 g_3 g_4}=a_{1234} + a_{1324} + a_{2314},
$$

\vspace{3mm}
$$
\sum_{\stackrel{(\{g_1, g_2\}, \{q_1, q_{2}\})}
{{}_{\{g_1, g_2, q_1, q_{2}\}=\{1, 2, 3, 4\}}}}
a_{g_1 g_2, q_1 q_{2}}=
$$

$$
=a_{12,34}+a_{13,24}+a_{14,23}
+a_{23,14}+a_{24,13}+a_{34,12},
$$

\vspace{3mm}
$$
\sum_{\stackrel{(\{g_1, g_2\}, \{q_1, q_{2}, q_3\})}
{{}_{\{g_1, g_2, q_1, q_{2}, q_3\}=\{1, 2, 3, 4, 5\}}}}
a_{g_1 g_2, q_1 q_{2}q_3}
=
$$

$$
=a_{12,345}+a_{13,245}+a_{14,235}
+a_{15,234}+a_{23,145}+a_{24,135}+
$$
$$
+a_{25,134}+a_{34,125}+a_{35,124}+a_{45,123},
$$

\vspace{4mm}
$$
\sum_{\stackrel{(\{\{g_1, g_2\}, \{g_3, g_{4}\}\}, \{q_1\})}
{{}_{\{g_1, g_2, g_3, g_{4}, q_1\}=\{1, 2, 3, 4, 5\}}}}
a_{g_1 g_2, g_3 g_{4},q_1}
=
$$

$$
=
a_{12,34,5}+a_{13,24,5}+a_{14,23,5}+
a_{12,35,4}+a_{13,25,4}+a_{15,23,4}+
$$
$$
+a_{12,54,3}+a_{15,24,3}+a_{14,25,3}+a_{15,34,2}+a_{13,54,2}+a_{14,53,2}+
$$
$$
+
a_{52,34,1}+a_{53,24,1}+a_{54,23,1}.
$$

\vspace{5mm}

Now we can write (\ref{tyyy}) as

\vspace{1mm}

$$
J[\psi^{(k)}]_{T,t}=
\hbox{\vtop{\offinterlineskip\halign{
\hfil#\hfil\cr
{\rm l.i.m.}\cr
$\stackrel{}{{}_{p_1,\ldots,p_k\to \infty}}$\cr
}} }
\sum\limits_{j_1=0}^{p_1}\ldots
\sum\limits_{j_k=0}^{p_k}
C_{j_k\ldots j_1}\Biggl(
\prod_{l=1}^k\zeta_{j_l}^{(i_l)}+\sum\limits_{r=1}^{[k/2]}
(-1)^r \times
\Biggr.
$$

\vspace{3mm}
\begin{equation}
\label{leto6000hh}
\times
\sum_{\stackrel{(\{\{g_1, g_2\}, \ldots, 
\{g_{2r-1}, g_{2r}\}\}, \{q_1, \ldots, q_{k-2r}\})}
{{}_{\{g_1, g_2, \ldots, 
g_{2r-1}, g_{2r}, q_1, \ldots, q_{k-2r}\}=\{1, 2, \ldots, k\}}}}
\prod\limits_{s=1}^r
{\bf 1}_{\{i_{g_{{}_{2s-1}}}=~i_{g_{{}_{2s}}}\ne 0\}}
\Biggl.{\bf 1}_{\{j_{g_{{}_{2s-1}}}=~j_{g_{{}_{2s}}}\}}
\prod_{l=1}^{k-2r}\zeta_{j_{q_l}}^{(i_{q_l})}\Biggr),
\end{equation}

\vspace{5mm}
\noindent
where $[x]$ is an integer part of a real number $x;$
another notations are the same as in Theorem {\bf 1}.

\vspace{2mm}

In particular, from (\ref{leto6000hh}) for $k=5$ we obtain

\vspace{3mm}

$$
J[\psi^{(5)}]_{T,t}=
\hbox{\vtop{\offinterlineskip\halign{
\hfil#\hfil\cr
{\rm l.i.m.}\cr
$\stackrel{}{{}_{p_1,\ldots,p_5\to \infty}}$\cr
}} }\sum_{j_1=0}^{p_1}\ldots\sum_{j_5=0}^{p_5}
C_{j_5\ldots j_1}\Biggl(
\prod_{l=1}^5\zeta_{j_l}^{(i_l)}-\Biggr.
$$

\vspace{2mm}
$$
-
\sum\limits_{\stackrel{(\{g_1, g_2\}, \{q_1, q_{2}, q_3\})}
{{}_{\{g_1, g_2, q_{1}, q_{2}, q_3\}=\{1, 2, 3, 4, 5\}}}}
{\bf 1}_{\{i_{g_{{}_{1}}}=~i_{g_{{}_{2}}}\ne 0\}}
{\bf 1}_{\{j_{g_{{}_{1}}}=~j_{g_{{}_{2}}}\}}
\prod_{l=1}^{3}\zeta_{j_{q_l}}^{(i_{q_l})}+
$$

\vspace{2mm}
$$
+
\sum_{\stackrel{(\{\{g_1, g_2\}, 
\{g_{3}, g_{4}\}\}, \{q_1\})}
{{}_{\{g_1, g_2, g_{3}, g_{4}, q_1\}=\{1, 2, 3, 4, 5\}}}}
{\bf 1}_{\{i_{g_{{}_{1}}}=~i_{g_{{}_{2}}}\ne 0\}}
{\bf 1}_{\{j_{g_{{}_{1}}}=~j_{g_{{}_{2}}}\}}
\Biggl.{\bf 1}_{\{i_{g_{{}_{3}}}=~i_{g_{{}_{4}}}\ne 0\}}
{\bf 1}_{\{j_{g_{{}_{3}}}=~j_{g_{{}_{4}}}\}}
\zeta_{j_{q_1}}^{(i_{q_1})}\Biggr).
$$

\vspace{7mm}
\noindent
The last equality obviously agrees with
(\ref{a5}).

Let us consider the generalization of Theorem 1 for the case
of an arbitrary complete orthonormal systems  
of functions in the space $L_2([t,T])$ 
and $\psi_1(\tau),\ldots,\psi_k(\tau)\in L_2([t, T]).$

\vspace{2mm}

{\bf Theorem~2}\ \cite{2018a} (Sect.~1.11), \cite{arxiv-1} (Sect.~15).
{\it Suppose that
$\psi_1(\tau),\ldots,\psi_k(\tau)\in L_2([t, T])$ and
$\{\phi_j(x)\}_{j=0}^{\infty}$ is an arbitrary complete orthonormal system  
of functions in the space $L_2([t,T]).$
Then the following expansion

\vspace{1mm}
$$
J[\psi^{(k)}]_{T,t}=
\hbox{\vtop{\offinterlineskip\halign{
\hfil#\hfil\cr
{\rm l.i.m.}\cr
$\stackrel{}{{}_{p_1,\ldots,p_k\to \infty}}$\cr
}} }
\sum\limits_{j_1=0}^{p_1}\ldots
\sum\limits_{j_k=0}^{p_k}
C_{j_k\ldots j_1}\Biggl(
\prod_{l=1}^k\zeta_{j_l}^{(i_l)}+\sum\limits_{r=1}^{[k/2]}
(-1)^r \times
\Biggr.
$$

\vspace{2mm}
\begin{equation}
\label{leto6000}
\times
\sum_{\stackrel{(\{\{g_1, g_2\}, \ldots, 
\{g_{2r-1}, g_{2r}\}\}, \{q_1, \ldots, q_{k-2r}\})}
{{}_{\{g_1, g_2, \ldots, 
g_{2r-1}, g_{2r}, q_1, \ldots, q_{k-2r}\}=\{1, 2, \ldots, k\}}}}
\prod\limits_{s=1}^r
{\bf 1}_{\{i_{g_{{}_{2s-1}}}=~i_{g_{{}_{2s}}}\ne 0\}}
\Biggl.{\bf 1}_{\{j_{g_{{}_{2s-1}}}=~j_{g_{{}_{2s}}}\}}
\prod_{l=1}^{k-2r}\zeta_{j_{q_l}}^{(i_{q_l})}\Biggr)
\end{equation}

\vspace{6mm}
\noindent
con\-verg\-ing in the mean-square sense is valid,
where $[x]$ is an integer part of a real number $x;$
another notations are the same as in Theorem~{\rm 1}.}

\vspace{2mm}

It should be noted that an analogue of Theorem 2 was considered 
in \cite{Rybakov1000}. 
Note that we use another notations 
\cite{2018a} (Sect.~1.11), \cite{arxiv-1} (Sect.~15)
in comparison with \cite{Rybakov1000}.
Moreover, the proof of an analogue of Theorem 2
from \cite{Rybakov1000} is somewhat different from the proof given in 
\cite{2018a} (Sect.~1.11), \cite{arxiv-1} (Sect.~15).

Note that, for the integrals $J[\psi^{(k)}]_{T,t}$ defined by
(\ref{ito}), 
the mean-square approximation error can be exactly
calculated and efficiently estimated.

Let $J[\psi^{(k)}]_{T,t}^{q}$ be the
expression on the right-hand side of (\ref{leto6000}) before passing to the limit 
for the case
$p_1=\ldots=p_k=q$, i.e.

\vspace{1mm}
$$
J[\psi^{(k)}]_{T,t}^q=
\sum\limits_{j_1,\ldots,j_k=0}^{q}
C_{j_k\ldots j_1}\Biggl(
\prod_{l=1}^k\zeta_{j_l}^{(i_l)}+\sum\limits_{r=1}^{[k/2]}
(-1)^r \times
\Biggr.
$$

\vspace{2mm}
\begin{equation}
\label{r1}
\times
\sum_{\stackrel{(\{\{g_1, g_2\}, \ldots, 
\{g_{2r-1}, g_{2r}\}\}, \{q_1, \ldots, q_{k-2r}\})}
{{}_{\{g_1, g_2, \ldots, 
g_{2r-1}, g_{2r}, q_1, \ldots, q_{k-2r}\}=\{1, 2, \ldots, k\}}}}
\prod\limits_{s=1}^r
{\bf 1}_{\{i_{g_{{}_{2s-1}}}=~i_{g_{{}_{2s}}}\ne 0\}}
\Biggl.{\bf 1}_{\{j_{g_{{}_{2s-1}}}=~j_{g_{{}_{2s}}}\}}
\prod_{l=1}^{k-2r}\zeta_{j_{q_l}}^{(i_{q_l})}\Biggr),
\end{equation}

\vspace{6mm}
\noindent
where $[x]$ is an integer part of a real number $x;$
another notations are the same as in Theorems~{\rm 1, 2}.

Let us denote

$$
{\sf M}\left\{\left(J[\psi^{(k)}]_{T,t}-
J[\psi^{(k)}]_{T,t}^{q}\right)^2\right\}\stackrel{{\rm def}}
{=}E_k^{q},
$$

\vspace{2mm}
$$
\int\limits_{[t,T]^k}
K^2(t_1,\ldots,t_k)dt_1\ldots dt_k
\stackrel{{\rm def}}{=}I_k.
$$

\vspace{4mm}

In \cite{2017-1}-\cite{2018aaa}, \cite{arxiv-1}, 
\cite{arxiv-2} it was shown that 

\vspace{-1mm}
\begin{equation}
\label{qq4}
E_k^{q}\le k!\left(I_k-\sum_{j_1,\ldots,j_k=0}^{q}C^2_{j_k\ldots j_1}\right)
\end{equation}

\vspace{2mm}
\noindent
for the following two cases:

\vspace{1mm}

1.\ $i_1,\ldots,i_k=1,\ldots,m$ and $T-t\in (0, +\infty)$,

2.\ $i_1,\ldots,i_k=0, 1,\ldots,m$ and  $T-t\in (0, 1)$.

\vspace{2mm}

The value $E_k^{q}$
can be calculated exactly.

\vspace{2mm}
              
{\bf Theorem 3} \cite{2018a} (Sect.~1.12), \cite{arxiv-2} (Sect.~6).
{\it Suppose that $\{\phi_j(x)\}_{j=0}^{\infty}$ 
is an arbitrary complete orthonormal system  
of functions in the space $L_2([t,T])$ and
$\psi_1(\tau),\ldots,\psi_k(\tau)\in L_2([t, T]),$  $i_1,\ldots, i_k=1,\ldots,m$.
Then

\begin{equation}
\label{tttr11}
E_k^p=I_k- \sum_{j_1,\ldots, j_k=0}^{p}
C_{j_k\ldots j_1}
{\sf M}\left\{J[\psi^{(k)}]_{T,t}
\sum\limits_{(j_1,\ldots,j_k)}
\int\limits_t^T \phi_{j_k}(t_k)
\ldots
\int\limits_t^{t_{2}}\phi_{j_{1}}(t_{1})
d{\bf f}_{t_1}^{(i_1)}\ldots
d{\bf f}_{t_k}^{(i_k)}\right\},
\end{equation}

\vspace{5mm}
\noindent
where $i_1,\ldots,i_k = 1,\ldots,m;$
the expression 

\vspace{-1mm}
$$
\sum\limits_{(j_1,\ldots,j_k)}
$$ 

\vspace{3mm}
\noindent
means the sum with respect to all
possible permutations 
$(j_1,\ldots,j_k)$. At the same time if 
$j_r$ swapped with $j_q$ in the permutation $(j_1,\ldots,j_k),$
then $i_r$ swapped with $i_q$ in the permutation
$(i_1,\ldots,i_k);$
another notations are the same as in Theorems {\rm 1, 2.}
}

\vspace{2mm}

Note that 

\vspace{-2mm}
$$
{\sf M}\left\{J[\psi^{(k)}]_{T,t}
\int\limits_t^T \phi_{j_k}(t_k)
\ldots
\int\limits_t^{t_{2}}\phi_{j_{1}}(t_{1})
d{\bf f}_{t_1}^{(i_1)}\ldots
d{\bf f}_{t_k}^{(i_k)}\right\}=C_{j_k\ldots j_1}.
$$

\vspace{4mm}

Therefore, for the case of pairwise 
different numbers $i_1,\ldots,i_k$ 
as well as for the case $i_1=\ldots=i_k$
from Theorem 3 it follows that
\cite{2018}-\cite{2018aaa}, \cite{17a}, 
\cite{arxiv-2}

\vspace{-2mm}
\begin{equation}
\label{qq1}
E_k^q= I_k- \sum_{j_1,\ldots,j_k=0}^{q}
C_{j_k\ldots j_1}^2,
\end{equation}

\vspace{1mm}
$$
E_k^q= I_k - \sum_{j_1,\ldots,j_k=0}^{q}
C_{j_k\ldots j_1}\Biggl(\sum\limits_{(j_1,\ldots,j_k)}
C_{j_k\ldots j_1}\Biggr),
$$

\vspace{4mm}
\noindent
where 
$$
\sum\limits_{(j_1,\ldots,j_k)}
$$ 

\vspace{3mm}
\noindent
is a sum with respect to all 
possible permutations
$(j_1,\ldots,j_k)$.

Consider some examples \cite{2018}-\cite{2018aaa}, \cite{17a}, 
\cite{arxiv-2} of application of Theorem 3 
$(i_1,i_2,i_3=1,\ldots,m)$

\vspace{1mm}
\begin{equation}
\label{qq2}
E_2^q     
=I_2
-\sum_{j_1,j_2=0}^q
C_{j_2j_1}^2-
\sum_{j_1,j_2=0}^q
C_{j_2j_1}C_{j_1j_2}\ \ \ (i_1=i_2),
\end{equation}

\vspace{2mm}
\begin{equation}
\label{qq3}
E_3^q=I_3
-\sum_{j_3,j_2,j_1=0}^q C_{j_3j_2j_1}^2-
\sum_{j_3,j_2,j_1=0}^q C_{j_3j_1j_2}C_{j_3j_2j_1}\ \ \ (i_1=i_2\ne i_3),
\end{equation}

\vspace{2mm}
\begin{equation}
\label{882}
E_3^q=I_3-
\sum_{j_3,j_2,j_1=0}^q C_{j_3j_2j_1}^2-
\sum_{j_3,j_2,j_1=0}^q C_{j_2j_3j_1}C_{j_3j_2j_1}\ \ \ (i_1\ne i_2=i_3),
\end{equation}

\vspace{2mm}
\begin{equation}
\label{883}
E_3^q=I_3
-\sum_{j_3,j_2,j_1=0}^q C_{j_3j_2j_1}^2-
\sum_{j_3,j_2,j_1=0}^q C_{j_3j_2j_1}C_{j_1j_2j_3}\ \ \ (i_1=i_3\ne i_2).
\end{equation}

\vspace{6mm}

The values $E_4^q$ and $E_5^q$ were calculated exaclty for all possible 
combinations of $i_1,\ldots,i_5=1,\ldots,m$ in 
\cite{2018}-\cite{2018aaa},  
\cite{arxiv-2}.

\vspace{5mm}

\section{Approximation of Specific Iterated Ito Stochastic Integrals
Based on Multiple Fourier--Legendre Series}

\vspace{5mm} 

Consider approximations of the iterated Ito stochastic integrals 
that appear in the numerical
scheme (\ref{4.45}) using Theorems 1, 2 for the case of 
complete orthonormal system of 
Legendre polynomials in the
space $L_2([\tau_p,\tau_{p+1}])$
($\tau_p=p\Delta,$ $N\Delta= T,$ 
$p=0,1,\ldots,N$) \cite{2006}
(also see \cite{2011-2}-\cite{arxiv-24}, \cite{Mikh-1}-\cite{Mikh-2})

\vspace{1mm}
\begin{equation}
\label{yyy1aaa}
I_{(0)\tau_{p+1},\tau_p}^{(i_1)}=\sqrt{\Delta}\zeta_0^{(i_1)},
\end{equation}

\vspace{2mm}
\begin{equation}
\label{qqqq1}
I_{(00)\tau_{p+1},\tau_p}^{(i_1 i_2)q}=
\frac{\Delta}{2}\left(\zeta_0^{(i_1)}\zeta_0^{(i_2)}+\sum_{i=1}^{q}
\frac{1}{\sqrt{4i^2-1}}\left(
\zeta_{i-1}^{(i_1)}\zeta_{i}^{(i_2)}-
\zeta_i^{(i_1)}\zeta_{i-1}^{(i_2)}\right) - {\bf 1}_{\{i_1=i_2\}}\right),
\end{equation}

\vspace{5mm}

\begin{equation}
\label{yyy2}
I_{(1)\tau_{p+1},\tau_p}^{(i_1)}=
-\frac{{\Delta}^{3/2}}{2}\left(\zeta_0^{(i_1)}+
\frac{1}{\sqrt{3}}\zeta_1^{(i_1)}\right),
\end{equation}

\vspace{5mm}

$$
I_{(000)\tau_{p+1},\tau_p}^{(i_1i_2i_3)q}
=\sum_{j_1,j_2,j_3=0}^{q}
C_{j_3j_2j_1}\Biggl(
\zeta_{j_1}^{(i_1)}\zeta_{j_2}^{(i_2)}\zeta_{j_3}^{(i_3)}
-{\bf 1}_{\{i_1=i_2\}}
{\bf 1}_{\{j_1=j_2\}}
\zeta_{j_3}^{(i_3)}-
\Biggr.
$$
\begin{equation}
\label{yyy3}
\Biggl.
-{\bf 1}_{\{i_2=i_3\}}
{\bf 1}_{\{j_2=j_3\}}
\zeta_{j_1}^{(i_1)}-
{\bf 1}_{\{i_1=i_3\}}
{\bf 1}_{\{j_1=j_3\}}
\zeta_{j_2}^{(i_2)}\Biggr),
\end{equation}

\vspace{5mm}

$$
I_{(0000)\tau_{p+1},\tau_p}^{(i_1 i_2 i_3 i_4)q}
=\sum_{j_1,j_2,j_3,j_4=0}^{q}
C_{j_4 j_3 j_2 j_1}\Biggl(
\prod_{l=1}^4\zeta_{j_l}^{(i_l)}
-\Biggr.
$$
$$
-
{\bf 1}_{\{i_1=i_2\}}
{\bf 1}_{\{j_1=j_2\}}
\zeta_{j_3}^{(i_3)}
\zeta_{j_4}^{(i_4)}
-
{\bf 1}_{\{i_1=i_3\}}
{\bf 1}_{\{j_1=j_3\}}
\zeta_{j_2}^{(i_2)}
\zeta_{j_4}^{(i_4)}-
$$
$$
-
{\bf 1}_{\{i_1=i_4\}}
{\bf 1}_{\{j_1=j_4\}}
\zeta_{j_2}^{(i_2)}
\zeta_{j_3}^{(i_3)}
-
{\bf 1}_{\{i_2=i_3\}}
{\bf 1}_{\{j_2=j_3\}}
\zeta_{j_1}^{(i_1)}
\zeta_{j_4}^{(i_4)}-
$$
$$
-
{\bf 1}_{\{i_2=i_4\}}
{\bf 1}_{\{j_2=j_4\}}
\zeta_{j_1}^{(i_1)}
\zeta_{j_3}^{(i_3)}
-
{\bf 1}_{\{i_3=i_4\}}
{\bf 1}_{\{j_3=j_4\}}
\zeta_{j_1}^{(i_1)}
\zeta_{j_2}^{(i_2)}+
$$
$$
+
{\bf 1}_{\{i_1=i_2\}}
{\bf 1}_{\{j_1=j_2\}}
{\bf 1}_{\{i_3=i_4\}}
{\bf 1}_{\{j_3=j_4\}}+
{\bf 1}_{\{i_1=i_3\}}
{\bf 1}_{\{j_1=j_3\}}
{\bf 1}_{\{i_2=i_4\}}
{\bf 1}_{\{j_2=j_4\}}+
$$
\begin{equation}
\label{yyy4}
+\Biggl.
{\bf 1}_{\{i_1=i_4\}}
{\bf 1}_{\{j_1=j_4\}}
{\bf 1}_{\{i_2=i_3\}}
{\bf 1}_{\{j_2=j_3\}}\Biggr),
\end{equation}

\vspace{6mm}

$$
I_{(01)\tau_{p+1},\tau_p}^{(i_1 i_2)q}=
-\frac{\Delta}{2}
I_{(00)\tau_{p+1},\tau_p}^{(i_1 i_2)q}
-\frac{{\Delta}^2}{4}\Biggl(
\frac{1}{\sqrt{3}}\zeta_0^{(i_1)}\zeta_1^{(i_2)}+\Biggr.
$$

\vspace{1mm}
\begin{equation}
\label{yyy5}
+\Biggl.\sum_{i=0}^{q}\Biggl(
\frac{(i+2)\zeta_i^{(i_1)}\zeta_{i+2}^{(i_2)}
-(i+1)\zeta_{i+2}^{(i_1)}\zeta_{i}^{(i_2)}}
{\sqrt{(2i+1)(2i+5)}(2i+3)}-
\frac{\zeta_i^{(i_1)}\zeta_{i}^{(i_2)}}{(2i-1)(2i+3)}\Biggr)\Biggr),
\end{equation}

\vspace{8mm}

$$
I_{(10)\tau_{p+1},\tau_p}^{(i_1 i_2)q}=
-\frac{\Delta}{2}I_{(00)\tau_{p+1},\tau_p}^{(i_1 i_2)q}
-\frac{\Delta^2}{4}\Biggl(
\frac{1}{\sqrt{3}}\zeta_0^{(i_2)}\zeta_1^{(i_1)}+\Biggr.
$$

\vspace{1mm}
\begin{equation}
\label{yyy6}
+\Biggl.\sum_{i=0}^{q}\Biggl(
\frac{(i+1)\zeta_{i+2}^{(i_2)}\zeta_{i}^{(i_1)}
-(i+2)\zeta_{i}^{(i_2)}\zeta_{i+2}^{(i_1)}}
{\sqrt{(2i+1)(2i+5)}(2i+3)}+
\frac{\zeta_i^{(i_1)}\zeta_{i}^{(i_2)}}{(2i-1)(2i+3)}\Biggr)\Biggr),
\end{equation}

\vspace{8mm}

\begin{equation}
\label{zzzz1}
{I}_{(2)\tau_{p+1},\tau_p}^{(i_1)}=
\frac{\Delta^{5/2}}{3}\left(
\zeta_0^{(i_1)}+\frac{\sqrt{3}}{2}\zeta_1^{(i_1)}+
\frac{1}{2\sqrt{5}}\zeta_2^{(i_1)}\right),
\end{equation}

\vspace{6mm}

$$
I_{(001)\tau_{p+1},\tau_p}^{(i_1i_2i_3)q}
=\sum_{j_1,j_2,j_3=0}^{q}
C_{j_3j_2j_1}^{001}\Biggl(
\zeta_{j_1}^{(i_1)}\zeta_{j_2}^{(i_2)}\zeta_{j_3}^{(i_3)}
-{\bf 1}_{\{i_1=i_2\}}
{\bf 1}_{\{j_1=j_2\}}
\zeta_{j_3}^{(i_3)}-
\Biggr.
$$
\begin{equation}
\label{yyy7}
\Biggl.
-{\bf 1}_{\{i_2=i_3\}}
{\bf 1}_{\{j_2=j_3\}}
\zeta_{j_1}^{(i_1)}-
{\bf 1}_{\{i_1=i_3\}}
{\bf 1}_{\{j_1=j_3\}}
\zeta_{j_2}^{(i_2)}\Biggr),
\end{equation}

\vspace{6mm}

$$
I_{(010)\tau_{p+1},\tau_p}^{(i_1i_2i_3)q}
=\sum_{j_1,j_2,j_3=0}^{q}
C_{j_3j_2j_1}^{010}\Biggl(
\zeta_{j_1}^{(i_1)}\zeta_{j_2}^{(i_2)}\zeta_{j_3}^{(i_3)}
-{\bf 1}_{\{i_1=i_2\}}
{\bf 1}_{\{j_1=j_2\}}
\zeta_{j_3}^{(i_3)}-
\Biggr.
$$
\begin{equation}
\label{yyy8}
\Biggl.
-{\bf 1}_{\{i_2=i_3\}}
{\bf 1}_{\{j_2=j_3\}}
\zeta_{j_1}^{(i_1)}-
{\bf 1}_{\{i_1=i_3\}}
{\bf 1}_{\{j_1=j_3\}}
\zeta_{j_2}^{(i_2)}\Biggr),
\end{equation}

\vspace{6mm}

$$
I_{(100)\tau_{p+1},\tau_p}^{(i_1i_2i_3)q}
=\sum_{j_1,j_2,j_3=0}^{q}
C_{j_3j_2j_1}^{100}\Biggl(
\zeta_{j_1}^{(i_1)}\zeta_{j_2}^{(i_2)}\zeta_{j_3}^{(i_3)}
-{\bf 1}_{\{i_1=i_2\}}
{\bf 1}_{\{j_1=j_2\}}
\zeta_{j_3}^{(i_3)}-
\Biggr.
$$
\begin{equation}
\label{yyy9}
\Biggl.
-{\bf 1}_{\{i_2=i_3\}}
{\bf 1}_{\{j_2=j_3\}}
\zeta_{j_1}^{(i_1)}-
{\bf 1}_{\{i_1=i_3\}}
{\bf 1}_{\{j_1=j_3\}}
\zeta_{j_2}^{(i_2)}\Biggr),
\end{equation}

\vspace{6mm}

$$
I_{(00000)\tau_{p+1},\tau_p}^{(i_1 i_2 i_3 i_4 i_5)q}
=\sum_{j_1,j_2,j_3,j_4,j_5=0}^q
C_{j_5 j_4 j_3 j_2 j_1}\Biggl(
\prod_{l=1}^5\zeta_{j_l}^{(i_l)}
-\Biggr.
$$
$$
-
{\bf 1}_{\{j_1=j_2\}}
{\bf 1}_{\{i_1=i_2\}}
\zeta_{j_3}^{(i_3)}
\zeta_{j_4}^{(i_4)}
\zeta_{j_5}^{(i_5)}-
{\bf 1}_{\{j_1=j_3\}}
{\bf 1}_{\{i_1=i_3\}}
\zeta_{j_2}^{(i_2)}
\zeta_{j_4}^{(i_4)}
\zeta_{j_5}^{(i_5)}-
$$
$$
-
{\bf 1}_{\{j_1=j_4\}}
{\bf 1}_{\{i_1=i_4\}}
\zeta_{j_2}^{(i_2)}
\zeta_{j_3}^{(i_3)}
\zeta_{j_5}^{(i_5)}-
{\bf 1}_{\{j_1=j_5\}}
{\bf 1}_{\{i_1=i_5\}}
\zeta_{j_2}^{(i_2)}
\zeta_{j_3}^{(i_3)}
\zeta_{j_4}^{(i_4)}-
$$
$$
-
{\bf 1}_{\{j_2=j_3\}}
{\bf 1}_{\{i_2=i_3\}}
\zeta_{j_1}^{(i_1)}
\zeta_{j_4}^{(i_4)}
\zeta_{j_5}^{(i_5)}-
{\bf 1}_{\{j_2=j_4\}}
{\bf 1}_{\{i_2=i_4\}}
\zeta_{j_1}^{(i_1)}
\zeta_{j_3}^{(i_3)}
\zeta_{j_5}^{(i_5)}-
$$
$$
-
{\bf 1}_{\{j_2=j_5\}}
{\bf 1}_{\{i_2=i_5\}}
\zeta_{j_1}^{(i_1)}
\zeta_{j_3}^{(i_3)}
\zeta_{j_4}^{(i_4)}-{\bf 1}_{\{j_3=j_4\}}
{\bf 1}_{\{i_3=i_4\}}
\zeta_{j_1}^{(i_1)}
\zeta_{j_2}^{(i_2)}
\zeta_{j_5}^{(i_5)}-
$$
$$
-
{\bf 1}_{\{j_3=j_5\}}
{\bf 1}_{\{i_3=i_5\}}
\zeta_{j_1}^{(i_1)}
\zeta_{j_2}^{(i_2)}
\zeta_{j_4}^{(i_4)}-{\bf 1}_{\{j_4=j_5\}}
{\bf 1}_{\{i_4=i_5\}}
\zeta_{j_1}^{(i_1)}
\zeta_{j_2}^{(i_2)}
\zeta_{j_3}^{(i_3)}+
$$
$$
+
{\bf 1}_{\{j_1=j_2\}}
{\bf 1}_{\{i_1=i_2\}}
{\bf 1}_{\{j_3=j_4\}}
{\bf 1}_{\{i_3=i_4\}}\zeta_{j_5}^{(i_5)}+
{\bf 1}_{\{j_1=j_2\}}
{\bf 1}_{\{i_1=i_2\}}
{\bf 1}_{\{j_3=j_5\}}
{\bf 1}_{\{i_3=i_5\}}\zeta_{j_4}^{(i_4)}+
$$
$$
+
{\bf 1}_{\{j_1=j_2\}}
{\bf 1}_{\{i_1=i_2\}}
{\bf 1}_{\{j_4=j_5\}}
{\bf 1}_{\{i_4=i_5\}}\zeta_{j_3}^{(i_3)}+
{\bf 1}_{\{j_1=j_3\}}
{\bf 1}_{\{i_1=i_3\}}
{\bf 1}_{\{j_2=j_4\}}
{\bf 1}_{\{i_2=i_4\}}\zeta_{j_5}^{(i_5)}+
$$
$$
+
{\bf 1}_{\{j_1=j_3\}}
{\bf 1}_{\{i_1=i_3\}}
{\bf 1}_{\{j_2=j_5\}}
{\bf 1}_{\{i_2=i_5\}}\zeta_{j_4}^{(i_4)}+
{\bf 1}_{\{j_1=j_3\}}
{\bf 1}_{\{i_1=i_3\}}
{\bf 1}_{\{j_4=j_5\}}
{\bf 1}_{\{i_4=i_5\}}\zeta_{j_2}^{(i_2)}+
$$
$$
+
{\bf 1}_{\{j_1=j_4\}}
{\bf 1}_{\{i_1=i_4\}}
{\bf 1}_{\{j_2=j_3\}}
{\bf 1}_{\{i_2=i_3\}}\zeta_{j_5}^{(i_5)}+
{\bf 1}_{\{j_1=j_4\}}
{\bf 1}_{\{i_1=i_4\}}
{\bf 1}_{\{j_2=j_5\}}
{\bf 1}_{\{i_2=i_5\}}\zeta_{j_3}^{(i_3)}+
$$
$$
+
{\bf 1}_{\{j_1=j_4\}}
{\bf 1}_{\{i_1=i_4\}}
{\bf 1}_{\{j_3=j_5\}}
{\bf 1}_{\{i_3=i_5\}}\zeta_{j_2}^{(i_2)}+
{\bf 1}_{\{j_1=j_5\}}
{\bf 1}_{\{i_1=i_5\}}
{\bf 1}_{\{j_2=j_3\}}
{\bf 1}_{\{i_2=i_3\}}\zeta_{j_4}^{(i_4)}+
$$
$$
+
{\bf 1}_{\{j_1=j_5\}}
{\bf 1}_{\{i_1=i_5\}}
{\bf 1}_{\{j_2=j_4\}}
{\bf 1}_{\{i_2=i_4\}}\zeta_{j_3}^{(i_3)}+
{\bf 1}_{\{j_1=j_5\}}
{\bf 1}_{\{i_1=i_5\}}
{\bf 1}_{\{j_3=j_4\}}
{\bf 1}_{\{i_3=i_4\}}\zeta_{j_2}^{(i_2)}+
$$
$$
+
{\bf 1}_{\{j_2=j_3\}}
{\bf 1}_{\{i_2=i_3\}}
{\bf 1}_{\{j_4=j_5\}}
{\bf 1}_{\{i_4=i_5\}}\zeta_{j_1}^{(i_1)}+
{\bf 1}_{\{j_2=j_4\}}
{\bf 1}_{\{i_2=i_4\}}
{\bf 1}_{\{j_3=j_5\}}
{\bf 1}_{\{i_3=i_5\}}\zeta_{j_1}^{(i_1)}+
$$
\begin{equation}
\label{yyy10}
+\Biggl.
{\bf 1}_{\{j_2=j_5\ne 0\}}
{\bf 1}_{\{i_2=i_5\}}
{\bf 1}_{\{j_3=j_4\ne 0\}}
{\bf 1}_{\{i_3=i_4\}}\zeta_{j_1}^{(i_1)}\Biggr),
\end{equation}

\vspace{7mm}
\noindent
where

\vspace{-1mm}
$$
C_{j_3j_2j_1}=\int\limits_{\tau_p}^{\tau_{p+1}}\phi_{j_3}(z)
\int\limits_{\tau_p}^{z}\phi_{j_2}(y)
\int\limits_{\tau_p}^{y}
\phi_{j_1}(x)dx dy dz=
$$

\begin{equation}
\label{hhh1}
=
\frac{\sqrt{(2j_1+1)(2j_2+1)(2j_3+1)}}{8}\Delta^{3/2}\bar
C_{j_3j_2j_1},
\end{equation}

\vspace{4mm}

$$
C_{j_4j_3j_2j_1}=\int\limits_{\tau_p}^{\tau_{p+1}}\phi_{j_4}(u)
\int\limits_{\tau_p}^{u}\phi_{j_3}(z)
\int\limits_{\tau_p}^{z}\phi_{j_2}(y)
\int\limits_{\tau_p}^{y}
\phi_{j_1}(x)dx dy dz du=
$$

\begin{equation}
\label{hhh2}
=\frac{\sqrt{(2j_1+1)(2j_2+1)(2j_3+1)(2j_4+1)}}{16}\Delta^{2}\bar
C_{j_4j_3j_2j_1},
\end{equation}

\vspace{4mm}

$$
C_{j_3j_2j_1}^{001}=\int\limits_{\tau_p}^{\tau_{p+1}}(\tau_p-z)\phi_{j_3}(z)
\int\limits_{\tau_p}^{z}\phi_{j_2}(y)
\int\limits_{\tau_p}^{y}
\phi_{j_1}(x)dx dy dz=
$$

\begin{equation}
\label{hhh3}
=
\frac{\sqrt{(2j_1+1)(2j_2+1)(2j_3+1)}}{16}\Delta^{5/2}\bar
C_{j_3j_2j_1}^{001},
\end{equation}

\vspace{4mm}

$$
C_{j_3j_2j_1}^{010}=\int\limits_{\tau_p}^{\tau_{p+1}}\phi_{j_3}(z)
\int\limits_{\tau_p}^{z}(\tau_p-y)\phi_{j_2}(y)
\int\limits_{\tau_p}^{y}
\phi_{j_1}(x)dx dy dz=
$$

\begin{equation}
\label{hhh4}
=
\frac{\sqrt{(2j_1+1)(2j_2+1)(2j_3+1)}}{16}\Delta^{5/2}\bar
C_{j_3j_2j_1}^{010},
\end{equation}

\vspace{4mm}

$$
C_{j_3j_2j_1}^{100}=\int\limits_{\tau_p}^{\tau_{p+1}}\phi_{j_3}(z)
\int\limits_{\tau_p}^{z}\phi_{j_2}(y)
\int\limits_{\tau_p}^{y}
(\tau_p-x)\phi_{j_1}(x)dx dy dz=
$$

\begin{equation}
\label{hhh5}
=
\frac{\sqrt{(2j_1+1)(2j_2+1)(2j_3+1)}}{16}\Delta^{5/2}\bar
C_{j_3j_2j_1}^{100},
\end{equation}

\vspace{4mm}

$$
C_{j_5j_4 j_3 j_2 j_1}=
\int\limits_{\tau_p}^{\tau_{p+1}}\phi_{j_5}(v)
\int\limits_{\tau_p}^v\phi_{j_4}(u)
\int\limits_{\tau_p}^{u}
\phi_{j_3}(z)
\int\limits_{\tau_p}^{z}\phi_{j_2}(y)\int\limits_{\tau_p}^{y}\phi_{j_1}(x)
dxdydzdudv=
$$

\begin{equation}
\label{hhh6}
=\frac{\sqrt{(2j_1+1)(2j_2+1)(2j_3+1)(2j_4+1)(2j_5+1)}}{32}\Delta^{5/2}\bar
C_{j_5j_4 j_3 j_2 j_1},
\end{equation}

\vspace{4mm}

\vspace{5mm}
\noindent
where
\begin{equation}
\label{jjj1}
\bar C_{j_3j_2j_1}=
\int\limits_{-1}^{1}P_{j_3}(z)
\int\limits_{-1}^{z}P_{j_2}(y)
\int\limits_{-1}^{y}
P_{j_1}(x)dx dy dz,
\end{equation}

\vspace{1mm}
\begin{equation}
\label{jjj2}
\bar C_{j_4j_3j_2j_1}=
\int\limits_{-1}^{1}P_{j_4}(u)
\int\limits_{-1}^{u}P_{j_3}(z)
\int\limits_{-1}^{z}P_{j_2}(y)
\int\limits_{-1}^{y}
P_{j_1}(x)dx dy dz,
\end{equation}

\vspace{1mm}

\begin{equation}
\label{jjj3}
\bar C_{j_3j_2j_1}^{100}=-
\int\limits_{-1}^{1}P_{j_3}(z)
\int\limits_{-1}^{z}P_{j_2}(y)
\int\limits_{-1}^{y}
P_{j_1}(x)(x+1)dx dy dz,
\end{equation}

\vspace{1mm}
\begin{equation}
\label{jjj4}
\bar C_{j_3j_2j_1}^{010}=-
\int\limits_{-1}^{1}P_{j_3}(z)
\int\limits_{-1}^{z}P_{j_2}(y)(y+1)
\int\limits_{-1}^{y}
P_{j_1}(x)dx dy dz,
\end{equation}

\vspace{1mm}

\begin{equation}
\label{jjj5}
\bar C_{j_3j_2j_1}^{001}=-
\int\limits_{-1}^{1}P_{j_3}(z)(z+1)
\int\limits_{-1}^{z}P_{j_2}(y)
\int\limits_{-1}^{y}
P_{j_1}(x)dx dy dz,
\end{equation}

\vspace{1mm}

\begin{equation}
\label{jjj6}
\bar C_{j_5j_4 j_3 j_2 j_1}=
\int\limits_{-1}^{1}P_{j_5}(v)
\int\limits_{-1}^{v}P_{j_4}(u)
\int\limits_{-1}^{u}P_{j_3}(z)
\int\limits_{-1}^{z}P_{j_2}(y)
\int\limits_{-1}^{y}
P_{j_1}(x)dx dy dz du dv,
\end{equation}

\vspace{6mm}
\noindent
where $P_i(x)$ $(i=0, 1, 2,\ldots)$ is the Legendre polynomial and

\vspace{2mm}

$$
\phi_i(x)=\sqrt{\frac{2i+1}{\Delta}}P_i\left(\left(x-\tau_p-\frac{\Delta}{2}\right)
\frac{2}{\Delta}\right),\ \ \ i=0, 1, 2,\ldots 
$$

\vspace{4mm}

Let us consider the exact relations and some estimates
for the mean-square errors of approximations of iterated Ito
stochastic integrals.

Using Theorem 3, we obtain
\cite{2017}-\cite{2013}, \cite{arxiv-3}, \cite{arxiv-24}

\vspace{1mm}
\begin{equation}
\label{xxx1}
{\sf M}\left\{\left(I_{(00)\tau_{p+1},\tau_p}^{(i_1 i_2)}-
I_{(00)\tau_{p+1},\tau_p}^{(i_1 i_2)q}
\right)^2\right\}
=\frac{\Delta^2}{2}\Biggl(\frac{1}{2}-\sum_{i=1}^q
\frac{1}{4i^2-1}\Biggr)\ \ \ (i_1\ne i_2),
\end{equation}

\vspace{7mm}

$$
{\sf M}\left\{\left(I_{(10)\tau_{p+1},\tau_p}^{(i_1 i_2)}-
I_{(10)\tau_{p+1},\tau_p}^{(i_1 i_2)q}
\right)^2\right\}=
{\sf M}\left\{\left(I_{(01)\tau_{p+1},\tau_p}^{(i_1 i_2)}-
I_{(01)\tau_{p+1},\tau_p}^{(i_1 i_2)q}\right)^2\right\}=
$$

\begin{equation}
\label{xxx2}
=\frac{\Delta^4}{16}\Biggl(\frac{5}{9}-
2\sum_{i=2}^q\frac{1}{4i^2-1}-
\sum_{i=1}^q
\frac{1}{(2i-1)^2(2i+3)^2}
-\sum_{i=0}^q\frac{(i+2)^2+(i+1)^2}{(2i+1)(2i+5)(2i+3)^2}
\Biggr)\ \ \ (i_1\ne i_2),
\end{equation}

\vspace{7mm}

$$
{\sf M}\left\{\left(I_{(10)\tau_{p+1},\tau_p}^{(i_1 i_1)}-
I_{(10)\tau_{p+1},\tau_p}^{(i_1 i_1)q}
\right)^2\right\}=
{\sf M}\left\{\left(I_{(01)\tau_{p+1},\tau_p}^{(i_1 i_1)}-
I_{(01)\tau_{p+1},\tau_p}^{(i_1 i_1)q}\right)^2\right\}=
$$

\begin{equation}
\label{xxx3}
=\frac{\Delta^4}{16}\Biggl(\frac{1}{9}-
\sum_{i=0}^{q}
\frac{1}{(2i+1)(2i+5)(2i+3)^2}
-2\sum_{i=1}^{q}
\frac{1}{(2i-1)^2(2i+3)^2}\Biggr).
\end{equation}

\vspace{9mm}

Applying (\ref{qq1}) and (\ref{qq2})--(\ref{883}), we get

\vspace{1mm}

$$
{\sf M}\left\{\left(
I_{(000)\tau_{p+1},\tau_p}^{(i_1i_2 i_3)}-
I_{(000)\tau_{p+1},\tau_p}^{(i_1i_2 i_3)q}\right)^2\right\}=
\frac{\Delta^{3}}{6}-\sum_{j_3,j_2,j_1=0}^{q}
C_{j_3j_2j_1}^2\ \ \ (i_1\ne i_2,\  i_1\ne i_3,\ i_2\ne i_3),
$$

\vspace{5mm}

$$
{\sf M}\left\{\left(
I_{(000)\tau_{p+1},\tau_p}^{(i_1i_2 i_3)}-
I_{(000)\tau_{p+1},\tau_p}^{(i_1i_2 i_3)q}\right)^2\right\}=
\frac{\Delta^{3}}{6}-\sum_{j_3,j_2,j_1=0}^{q}
C_{j_3j_2j_1}^2
-\sum_{j_3,j_2,j_1=0}^{q}
C_{j_2j_3j_1}C_{j_3j_2j_1}\ \ \ (i_1\ne i_2=i_3),
$$

\vspace{5mm}

$$
{\sf M}\left\{\left(
I_{(000)\tau_{p+1},\tau_p}^{(i_1i_2 i_3)}-
I_{(000)\tau_{p+1},\tau_p}^{(i_1i_2 i_3)q}\right)^2\right\}=
\frac{\Delta^{3}}{6}-\sum_{j_3,j_2,j_1=0}^{q}
C_{j_3j_2j_1}^2
-\sum_{j_3,j_2,j_1=0}^{q}
C_{j_3j_2j_1}C_{j_1j_2j_3}\ \ \ (i_1=i_3\ne i_2),
$$

\vspace{5mm}

$$
{\sf M}\left\{\left(
I_{(000)\tau_{p+1},\tau_p}^{(i_1i_2 i_3)}-
I_{(000)\tau_{p+1},\tau_p}^{(i_1i_2 i_3)q}\right)^2\right\}=
\frac{\Delta^{3}}{6}-\sum_{j_3,j_2,j_1=0}^{q}
C_{j_3j_2j_1}^2
-\sum_{j_3,j_2,j_1=0}^{q}
C_{j_3j_1j_2}C_{j_3j_2j_1}\ \ \ (i_1=i_2\ne i_3).
$$

\vspace{8mm}

At the same time using the estimate (\ref{qq4}) 
for $i_1,\ldots,i_5=1,\ldots,m$, 
we have

\vspace{2mm}

$$
{\sf M}\left\{\left(
I_{(01)\tau_{p+1},\tau_p}^{(i_1i_2)}-
I_{(01)\tau_{p+1},\tau_p}^{(i_1i_2)q}\right)^2\right\}\le
2\Biggl(\frac{\Delta^{4}}{4}-\sum_{j_1,j_2=0}^{q}
\left(C_{j_2j_1}^{01}\right)^2\Biggr),
$$

\vspace{3mm}
$$
{\sf M}\left\{\left(
I_{(10)\tau_{p+1},\tau_p}^{(i_1i_2)}-
I_{(10)\tau_{p+1},\tau_p}^{(i_1i_2)q}\right)^2\right\}\le
2\Biggl(\frac{\Delta^{4}}{12}-\sum_{j_1,j_2=0}^{q}
\left(C_{j_2j_1}^{10}\right)^2\Biggr),
$$

\vspace{3mm}

\begin{equation}
\label{xxx4}
{\sf M}\left\{\left(
I_{(000)\tau_{p+1},\tau_p}^{(i_1i_2 i_3)}-
I_{(000)\tau_{p+1},\tau_p}^{(i_1i_2 i_3)q}\right)^2\right\}\le
6\Biggl(\frac{\Delta^{3}}{6}-\sum_{j_3,j_2,j_1=0}^{q}
C_{j_3j_2j_1}^2\Biggr),
\end{equation}

\vspace{2mm}

\begin{equation}
\label{xxx5}
{\sf M}\left\{\left(
I_{(0000)\tau_{p+1},\tau_p}^{(i_1i_2 i_3 i_4)}-
I_{(0000)\tau_{p+1},\tau_p}^{(i_1i_2 i_3 i_4)q}\right)^2\right\}\le
24\Biggl(\frac{\Delta^{4}}{24}-\sum_{j_1,j_2,j_3,j_4=0}^{q}
C_{j_4j_3j_2j_1}^2\Biggr),
\end{equation}

\vspace{2mm}

\begin{equation}
\label{xxx6}
{\sf M}\left\{\left(
I_{(100)\tau_{p+1},\tau_p}^{(i_1i_2 i_3)}-
I_{(100)\tau_{p+1},\tau_p}^{(i_1i_2 i_3)q}\right)^2\right\}\le
6\Biggl(\frac{\Delta^{5}}{60}-\sum_{j_1,j_2,j_3=0}^{q}
\left(C_{j_3j_2j_1}^{100}\right)^2\Biggr),
\end{equation}

\vspace{2mm}

\begin{equation}
\label{xxx7}
{\sf M}\left\{\left(
I_{(010)\tau_{p+1},\tau_p}^{(i_1i_2 i_3)}-
I_{(010)\tau_{p+1},\tau_p}^{(i_1i_2 i_3)q}\right)^2\right\}\le
6\Biggl(\frac{\Delta^{5}}{20}-\sum_{j_1,j_2,j_3=0}^{q}
\left(C_{j_3j_2j_1}^{010}\right)^2\Biggr),
\end{equation}
                          
\vspace{2mm}

\begin{equation}
\label{xxx8}
{\sf M}\left\{\left(
I_{(001)\tau_{p+1},\tau_p}^{(i_1i_2 i_3)}-
I_{(001)\tau_{p+1},\tau_p}^{(i_1i_2 i_3)q}\right)^2\right\}\le
6\Biggl(\frac{\Delta^5}{10}-\sum_{j_1,j_2,j_3=0}^{q}
\left(C_{j_3j_2j_1}^{001}\right)^2\Biggr),
\end{equation}

\vspace{2mm}

\begin{equation}
\label{xxx9}
{\sf M}\left\{\left(
I_{(00000)\tau_{p+1},\tau_p}^{(i_1 i_2 i_3 i_4 i_5)}-
I_{(00000)\tau_{p+1},\tau_p}^{(i_1 i_2 i_3 i_4 i_5)q}\right)^2\right\}\le
120\left(\frac{\Delta^{5}}{120}-\sum_{j_1,j_2,j_3,j_4,j_5=0}^{q}
C_{j_5 i_4 i_3 i_2 j_1}^2\right).
\end{equation}

\vspace{5mm}

The Fourier--Legendre coefficients

$$
\bar C_{j_3j_2j_1},\ \ \ \bar C_{j_4j_3j_2j_1},\ \ \
\bar C_{j_3j_2j_1}^{001},\ \ \ \bar C_{j_3j_2j_1}^{010},\ \ \
\bar C_{j_3j_2j_1}^{100},\ \ \
\bar C_{j_5j_4j_3j_2j_1}
$$ 

\vspace{4mm}
\noindent
can be calculated exactly using computer algebra 
systems like
Derive. The exact values of these Fourier--Legendre coefficients  
were presented in
tabular form in the monographs \cite{2006}-\cite{2013}. 
Note that the mendioned Fourier--Legendre coefficients 
do not depend
on the integration step $\tau_{p+1}-\tau_p$ of the numerical method, 
which can 
be variable.

Recently, 
the database with 270,000 exactly
calculated Fourier--Legendre coefficients was described \cite{Mikh-1}, \cite{Kuz-Kuz}.
This database was used in the software package,
which is written in the Python programming language
for the implementation of explicit one-step strong numerical schemes 
with orders 0.5, 1.0, 1.5, 2.0, 2.5, and 3.0 
of convergence for Ito stochastic differential equations
\cite{Mikh-1}, \cite{Kuz-Kuz}.
The optimization of the mean-square approximation 
procedures for iterated Ito stochastic integrals
from these numerical schemes can be found in \cite{Mikh-2}.

Note that in \cite{KlPl2}-\cite{Mi3} (also see \cite{rr})
the iterated stochastic integrals were 
approximated using the trigonometric
Fourier expansion of the multidimensional Brownian bridge process. 
It is important to pay attention 
that the number $q$ must be the same for all
approximations of iterated stochastic integrals from the
considered collection in the approcah from 
\cite{KlPl2}-\cite{Mi3}. At the same time
the numbers $q$ can be 
selected individually for different stochastic integrals
from the
considred collection in the method based on Theorems 1--3
(see Sect.~5.3, 6.2 from \cite{2018a} for details).

On the basis of 
the presented 
expansions (\ref{yyy1aaa})--(\ref{yyy10}) of 
iterated Ito stochastic integrals we 
can see that increasing of multiplicities of these integrals 
or degree indexes of their weight functions 
leads
to increasing 
of smallness orders with respect to $\Delta$ in the mean-square sense 
for iterated stochastic integrals. This leads to a sharp decrease  
of member 
quantities (the numbers $q$)
in expansions of iterated Ito stochastic 
integrals, which are required for achieving the acceptable accuracy
of approximation.
Generally speaking, the minimum values $q$ that guarantee the 
fulfillment of the condition 
(\ref{4.3})
for each approximation (see (\ref{yyy1aaa})--(\ref{yyy10}))
are different and abruptly decreasing with the growth of 
smallness order (with respect to $\Delta$) of
the approximations of iterated stochastic integrals.

The detailed comparison of the method from 
\cite{KlPl2}-\cite{Mi3} with the method based on Theorems 1--3
can be found in \cite{2018a}-\cite{2018aaa}
(Chapters 2, 5, 6), \cite{29a}, \cite{301a}.

\vspace{5mm}

\section{Approximation of Iterated Ito and Stratonovich Stochastic
Integrals. Combined Approach}

\vspace{5mm}

In contrast to the iterated Ito stochastic integrals (\ref{ito}), 
the iterated Stratonovich stochastic integrals (\ref{str})
have simpler expansions (see Theorems 4--10 below) 
than (\ref{leto6000}). However, the calculation (or estimation) 
of the mean-square approximation
error for the latter is a more difficult problem than 
for the former. Below in this section, we will study this issue in
detail.

As it turned out, Theorems 1, 2 can be adapted for the iterated
Stratonovich stochastic integrals (\ref{str}) at least
for multiplicities 1 to 6. 
Expansions of these iterated Stratonovich 
stochastic integrals turned out
much simpler than the appropriate expansions
of the iterated Ito stochastic integrals (\ref{ito}) from Theorems 1, 2.
Applying this feature and standard relations between
iterated Ito and Stratonovich stochastic integrals, we will get 
simpler expansions for the 
iterated Ito stochastic integrals (\ref{ito}) than the
expansions from the previous section. However, as was mentioned
above, the estimation
of the mean-square approximation
error for the expansions from this section is a nontrivial problem.

Let us first present some old results on expansion of the iterated
Stratonovich stochastic integrals (\ref{str}) of
multiplicities 2 to 4 (Theorems 4--6 below).

\vspace{2mm}

{\bf Theorem 4}\ \cite{2011-2}-\cite{2018aaa}, 
\cite{2010-2}-\cite{2013}, \cite{30a}, \cite{300a},
\cite{400a}, \cite{271a},
\cite{arxiv-5}, \cite{arxiv-8}, \cite{arxiv-23}.\
{\it Assume that the following conditions are fulfilled{\rm :}

{\rm 1}. The function $\psi_2(\tau)$ is continuously 
differentiable at the interval $[t, T]$ and the
function $\psi_1(\tau)$ is twice continuously 
differentiable at the interval $[t, T]$.

{\rm 2}. $\{\phi_j(x)\}_{j=0}^{\infty}$ is a complete orthonormal system 
of Legendre polynomials or tri\-go\-no\-met\-ric functions
in the space $L_2([t, T]).$

Then, the iterated Stratonovich 
stochastic integral of the second multiplicity

$$
{\int\limits_t^{*}}^T\psi_2(t_2)
{\int\limits_t^{*}}^{t_2}\psi_1(t_1)d{\bf f}_{t_1}^{(i_1)}
d{\bf f}_{t_2}^{(i_2)}\ \ \ (i_1, i_2=1,\ldots,m)
$$

\vspace{3mm}
\noindent
is expanded into the 
double series

$$
{\int\limits_t^{*}}^T\psi_2(t_2)
{\int\limits_t^{*}}^{t_2}\psi_1(t_1)d{\bf f}_{t_1}^{(i_1)}
d{\bf f}_{t_2}^{(i_2)}=
\hbox{\vtop{\offinterlineskip\halign{
\hfil#\hfil\cr
{\rm l.i.m.}\cr
$\stackrel{}{{}_{p_1,p_2\to \infty}}$\cr
}} }\sum_{j_1=0}^{p_1}\sum_{j_2=0}^{p_2}
C_{j_2j_1}\zeta_{j_1}^{(i_1)}\zeta_{j_2}^{(i_2)}
$$

\vspace{3mm}
\noindent
converging in the mean-square sense, where 

\vspace{1mm}
$$
C_{j_2 j_1}=\int\limits_t^T
\psi_2(s_2)\phi_{j_2}(s_2)
\int\limits_t^{s_2}\psi_1(s_1)\phi_{j_1}(s_1)ds_1ds_2;
$$

\vspace{4mm}
\noindent
another notations are the same as in Theorems {\rm 1, 2}.}

\vspace{2mm}

{\bf Theorem 5}\ \cite{2011-2}-\cite{2018aaa}, 
\cite{2010-2}-\cite{2013}, \cite{271a},
\cite{arxiv-5}, \cite{arxiv-7}.\
{\it Assume that
$\{\phi_j(x)\}_{j=0}^{\infty}$ is a complete orthonormal
system of Legendre polynomials or trigonomertic functions
in the space $L_2([t, T])$. Furthermore,
the function $\psi_2(\tau)$ is continuously
differentiable at the interval $[t, T]$ and
the functions $\psi_1(\tau),$ $\psi_3(\tau)$ are twice continuously
differentiable at the interval $[t, T]$.
Then, for the iterated Stratonovich stochastic integral of third multiplicity

$$
{\int\limits_t^{*}}^T\psi_3(t_3)
{\int\limits_t^{*}}^{t_3}\psi_2(t_2)
{\int\limits_t^{*}}^{t_2}\psi_1(t_1)
d{\bf f}_{t_1}^{(i_1)}
d{\bf f}_{t_2}^{(i_2)}d{\bf f}_{t_3}^{(i_3)}\ \ \ (i_1, i_2, i_3=1,\ldots,m)
$$

\vspace{3mm}
\noindent
the following expansion

\begin{equation}
\label{feto19000a}
{\int\limits_t^{*}}^T\psi_3(t_3)
{\int\limits_t^{*}}^{t_3}\psi_2(t_2)
{\int\limits_t^{*}}^{t_2}\psi_1(t_1)
d{\bf f}_{t_1}^{(i_1)}
d{\bf f}_{t_2}^{(i_2)}d{\bf f}_{t_3}^{(i_3)}\ 
=
\hbox{\vtop{\offinterlineskip\halign{
\hfil#\hfil\cr
{\rm l.i.m.}\cr
$\stackrel{}{{}_{q\to \infty}}$\cr
}} }
\sum\limits_{j_1, j_2, j_3=0}^{q}
C_{j_3 j_2 j_1}\zeta_{j_1}^{(i_1)}\zeta_{j_2}^{(i_2)}\zeta_{j_3}^{(i_3)}
\end{equation}

\vspace{3mm}
\noindent
converging in the mean-square sense 
is valid, where

\vspace{1mm}
$$
C_{j_3 j_2 j_1}=\int\limits_t^T\psi_3(t_3)\phi_{j_3}(t_3)
\int\limits_t^{t_3}\psi_2(t_2)\phi_{j_2}(t_2)
\int\limits_t^{t_2}\psi_1(t_1)\phi_{j_1}(t_1)dt_1dt_2dt_3;
$$

\vspace{4mm}
\noindent
another notations are the same as in Theorems {\rm 1, 2}.}

\vspace{2mm}

{\bf Theorem 6}\ \cite{2011-2}-\cite{2018aaa}, 
\cite{2010-2}-\cite{2013}, \cite{271a},
\cite{arxiv-5}, \cite{arxiv-6}.\ 
{\it Suppose that 
$\{\phi_j(x)\}_{j=0}^{\infty}$ is a complete orthonormal system of 
Legendre polynomials or trigonometric functions in $L_2([t, T]).$
Then, for the iterated Stra\-to\-no\-vich stochastic integral
of multiplicity {\rm 4} 

\vspace{1mm}
$$
I_{T,t}^{*(i_1 i_2 i_3 i_4)}=
{\int\limits_t^{*}}^T
{\int\limits_t^{*}}^{t_4}
{\int\limits_t^{*}}^{t_3}
{\int\limits_t^{*}}^{t_2}
d{\bf w}_{t_1}^{(i_1)}
d{\bf w}_{t_2}^{(i_2)}d{\bf w}_{t_3}^{(i_3)}d{\bf w}_{t_4}^{(i_4)}
$$

\vspace{4mm}
\noindent
the following 
expansion

$$
I_{T,t}^{*(i_1 i_2 i_3 i_4)}=
\hbox{\vtop{\offinterlineskip\halign{
\hfil#\hfil\cr
{\rm l.i.m.}\cr
$\stackrel{}{{}_{q\to \infty}}$\cr
}} }
\sum\limits_{j_1, j_2, j_3, j_4=0}^{q}
C_{j_4 j_3 j_2 j_1}\zeta_{j_1}^{(i_1)}\zeta_{j_2}^{(i_2)}\zeta_{j_3}^{(i_3)}
\zeta_{j_4}^{(i_4)}
$$

\vspace{5mm}
\noindent
converging in the mean-square sense 
is valid, where $i_1, i_2, i_3, i_4=0, 1,\ldots,m,$

\vspace{1mm}
$$
C_{j_4 j_3 j_2 j_1}=\int\limits_t^T\phi_{j_4}(t_4)\int\limits_t^{t_4}
\phi_{j_3}(t_3)
\int\limits_t^{t_3}\phi_{j_2}(t_2)\int\limits_t^{t_2}\phi_{j_1}(t_1)
dt_1dt_2dt_3dt_4,
$$

\vspace{4mm}
\noindent
${\bf w}_{\tau}^{(i)}={\bf f}_{\tau}^{(i)}$ $(i=1,\ldots,m)$ are independent 
standard Wiener processes and 
${\bf w}_{\tau}^{(0)}=\tau.$}

\vspace{2mm}

Recently, a new approach to the expansion and mean-square 
approximation of iterated Stratonovich stochastic integrals has been obtained
\cite{2018a} (Sect.~2.10--2.16), 
\cite{arxiv-4} (Sect.~7--13),
\cite{arxiv-5} (Sect.~13--19), 
\cite{arxiv-6} (Sect.~5--11), \cite{new-art-1-xxy}
(Sect.~4--9), \cite{new-art-1xxys}.
Let us formulate four theorems that were obtained using this approach.

\vspace{2mm}

{\bf Theorem 7}\ \cite{2018a}, \cite{arxiv-4}, \cite{arxiv-5}, \cite{arxiv-6},
\cite{new-art-1-xxy}.\
{\it Suppose 
that $\{\phi_j(x)\}_{j=0}^{\infty}$ is a complete orthonormal system of 
Legendre polynomials or trigonometric functions in the space $L_2([t, T]).$
Furthermore, let $\psi_1(\tau), \psi_2(\tau),$ $\psi_3(\tau)$ are continuously dif\-ferentiable 
nonrandom functions on $[t, T].$ 
Then, for the 
iterated Stra\-to\-no\-vich stochastic integral of third multiplicity

$$
J^{*}[\psi^{(3)}]_{T,t}={\int\limits_t^{*}}^T\psi_3(t_3)
{\int\limits_t^{*}}^{t_3}\psi_2(t_2)
{\int\limits_t^{*}}^{t_2}\psi_1(t_1)
d{\bf w}_{t_1}^{(i_1)}
d{\bf w}_{t_2}^{(i_2)}d{\bf w}_{t_3}^{(i_3)}\ \ \ (i_1,i_2,i_3=0,1,\ldots,m)
$$

\vspace{4mm}
\noindent
the following 
relations

\vspace{-1mm}
\begin{equation}
\label{fin1}
J^{*}[\psi^{(3)}]_{T,t}
=\hbox{\vtop{\offinterlineskip\halign{
\hfil#\hfil\cr
{\rm l.i.m.}\cr
$\stackrel{}{{}_{p\to \infty}}$\cr
}} }
\sum\limits_{j_1, j_2, j_3=0}^{p}
C_{j_3 j_2 j_1}\zeta_{j_1}^{(i_1)}\zeta_{j_2}^{(i_2)}\zeta_{j_3}^{(i_3)},
\end{equation}

\vspace{3mm}
\begin{equation}
\label{fin2}
{\sf M}\left\{\left(
J^{*}[\psi^{(3)}]_{T,t}-
\sum\limits_{j_1, j_2, j_3=0}^{p}
C_{j_3 j_2 j_1}\zeta_{j_1}^{(i_1)}\zeta_{j_2}^{(i_2)}\zeta_{j_3}^{(i_3)}\right)^2\right\}
\le \frac{C}{p}
\end{equation}

\vspace{5mm}
\noindent
are fulfilled, where $i_1, i_2, i_3=0,1,\ldots,m$ in {\rm (\ref{fin1})} and 
$i_1, i_2, i_3=1,\ldots,m$ in {\rm (\ref{fin2})},
constant $C$ is independent of $p,$

\vspace{1mm}
$$
C_{j_3 j_2 j_1}=\int\limits_t^T\psi_3(t_3)\phi_{j_3}(t_3)
\int\limits_t^{t_3}\psi_2(t_2)\phi_{j_2}(t_2)
\int\limits_t^{t_2}\psi_1(t_1)\phi_{j_1}(t_1)dt_1dt_2dt_3
$$

\vspace{4mm}
\noindent
and
$$
\zeta_{j}^{(i)}=
\int\limits_t^T \phi_{j}(\tau) d{\bf f}_{\tau}^{(i)}
$$ 

\vspace{2mm}
\noindent
are independent standard Gaussian random variables for various 
$i$ or $j$ {\rm (}in the case when $i\ne 0${\rm );} 
another notations are the same as in Theorems~{\rm 1, 2}.}

\vspace{2mm}

{\bf Theorem 8}\ \cite{2018a}, \cite{arxiv-4}, \cite{arxiv-5}, \cite{arxiv-6},
\cite{new-art-1-xxy}.\ 
{\it Let
$\{\phi_j(x)\}_{j=0}^{\infty}$ be a complete orthonormal system of 
Legendre polynomials or trigonometric functions in the space $L_2([t, T]).$
Furthermore, let $\psi_1(\tau), \ldots, \psi_4(\tau)$ be continuously dif\-ferentiable 
nonrandom functions on $[t, T].$ 
Then, for the 
iterated Stra\-to\-no\-vich stochastic integral of fourth multiplicity

\begin{equation}
\label{fin0}
J^{*}[\psi^{(4)}]_{T,t}={\int\limits_t^{*}}^T\psi_4(t_4)
{\int\limits_t^{*}}^{t_4}\psi_3(t_3)
{\int\limits_t^{*}}^{t_3}\psi_2(t_2)
{\int\limits_t^{*}}^{t_2}\psi_1(t_1)
d{\bf w}_{t_1}^{(i_1)}
d{\bf w}_{t_2}^{(i_2)}d{\bf w}_{t_3}^{(i_3)}d{\bf w}_{t_4}^{(i_4)}
\end{equation}

\vspace{4mm}
\noindent
the following 
relations

\begin{equation}
\label{fin3}
J^{*}[\psi^{(4)}]_{T,t}
=\hbox{\vtop{\offinterlineskip\halign{
\hfil#\hfil\cr
{\rm l.i.m.}\cr
$\stackrel{}{{}_{p\to \infty}}$\cr
}} }
\sum\limits_{j_1, j_2, j_3,j_4=0}^{p}
C_{j_4j_3 j_2 j_1}\zeta_{j_1}^{(i_1)}\zeta_{j_2}^{(i_2)}\zeta_{j_3}^{(i_3)}\zeta_{j_4}^{(i_4)},
\end{equation}

\vspace{3mm}

\begin{equation}
\label{fin4}
{\sf M}\left\{\left(
J^{*}[\psi^{(4)}]_{T,t}-
\sum\limits_{j_1, j_2, j_3, j_4=0}^{p}
C_{j_4 j_3 j_2 j_1}\zeta_{j_1}^{(i_1)}\zeta_{j_2}^{(i_2)}\zeta_{j_3}^{(i_3)}
\zeta_{j_4}^{(i_4)}
\right)^2\right\}
\le \frac{C}{p^{1-\varepsilon}}
\end{equation}

\vspace{5mm}
\noindent
are fulfilled, where $i_1, \ldots , i_4=0,1,\ldots,m$ in {\rm (\ref{fin0}),} {\rm (\ref{fin3})} 
and $i_1, \ldots, i_4=1,\ldots,m$ in {\rm (\ref{fin4}),}
constant $C$ does not depend on $p,$
$\varepsilon$ is an arbitrary
small positive real number 
for the case of complete orthonormal system of 
Legendre polynomials in the space $L_2([t, T])$
and $\varepsilon=0$ for the case of
complete orthonormal system of 
trigonometric functions in the space $L_2([t, T]),$

\vspace{1mm}
$$
C_{j_4 j_3 j_2 j_1}=
$$

$$
=
\int\limits_t^T\psi_4(t_4)\phi_{j_4}(t_4)
\int\limits_t^{t_4}\psi_3(t_3)\phi_{j_3}(t_3)
\int\limits_t^{t_3}\psi_2(t_2)\phi_{j_2}(t_2)
\int\limits_t^{t_2}\psi_1(t_1)\phi_{j_1}(t_1)dt_1dt_2dt_3dt_4;
$$

\vspace{4mm}
\noindent
another notations are the same as in Theorem~{\rm 7}.}

\vspace{2mm}

{\bf Theorem 9}\ \cite{2018a}, \cite{arxiv-4}, \cite{arxiv-5}, \cite{arxiv-6},
\cite{new-art-1-xxy}.\
{\it Assume 
that $\{\phi_j(x)\}_{j=0}^{\infty}$ is a complete orthonormal system of 
Legendre polynomials or trigonometric functions in the space $L_2([t, T])$
and $\psi_1(\tau), \ldots, \psi_5(\tau)$ are continuously dif\-ferentiable 
nonrandom functions on $[t, T].$ 
Then, for the 
iterated Stra\-to\-no\-vich stochastic integral of fifth multiplicity

\vspace{1mm}
\begin{equation}
\label{fin7}
J^{*}[\psi^{(5)}]_{T,t}={\int\limits_t^{*}}^T\psi_5(t_5)
\ldots
{\int\limits_t^{*}}^{t_2}\psi_1(t_1)
d{\bf w}_{t_1}^{(i_1)}
\ldots d{\bf w}_{t_5}^{(i_5)}
\end{equation}

\vspace{5mm}
\noindent
the following 
relations

\begin{equation}
\label{fin8}
J^{*}[\psi^{(5)}]_{T,t}
=\hbox{\vtop{\offinterlineskip\halign{
\hfil#\hfil\cr
{\rm l.i.m.}\cr
$\stackrel{}{{}_{p\to \infty}}$\cr
}} }
\sum\limits_{j_1,\ldots,j_5=0}^{p}
C_{j_5 \ldots j_1}\zeta_{j_1}^{(i_1)}\ldots \zeta_{j_5}^{(i_5)},
\end{equation}

\vspace{3mm}

\begin{equation}
\label{fin9}
{\sf M}\left\{\left(
J^{*}[\psi^{(5)}]_{T,t}-
\sum\limits_{j_1, \ldots, j_5=0}^{p}
C_{j_5 \ldots j_1}\zeta_{j_1}^{(i_1)}\ldots
\zeta_{j_5}^{(i_5)}
\right)^2\right\}
\le \frac{C}{p^{1-\varepsilon}}
\end{equation}

\vspace{5mm}
\noindent
are fulfilled, where $i_1, \ldots , i_5=0,1,\ldots,m$ in {\rm (\ref{fin7}),} {\rm (\ref{fin8})} 
and $i_1, \ldots, i_5=1,\ldots,m$ in {\rm (\ref{fin9}),}
constant $C$ is independent of $p,$
$\varepsilon$ is an arbitrary
small positive real number 
for the case of complete orthonormal system of 
Legendre polynomials in the space $L_2([t, T])$
and $\varepsilon=0$ for the case of
complete orthonormal system of 
trigonometric functions in the space $L_2([t, T]),$

\vspace{1mm}
$$
C_{j_5 \ldots j_1}=
\int\limits_t^T\psi_5(t_5)\phi_{j_5}(t_5)\ldots
\int\limits_t^{t_2}\psi_1(t_1)\phi_{j_1}(t_1)dt_1\ldots dt_5;
$$

\vspace{4mm}
\noindent
another notations are the same as in Theorems~{\rm 7, 8}.}

\vspace{2mm}

{\bf Theorem 10}\ \cite{2018a}, \cite{arxiv-4}, \cite{arxiv-5}, \cite{arxiv-6},
\cite{new-art-1xxys}.\
{\it Suppose that 
$\{\phi_j(x)\}_{j=0}^{\infty}$ is a complete orthonormal system of 
Legendre polynomials or trigonometric functions in the space $L_2([t, T]).$
Then, for the 
iterated Stratonovich stochastic integral of sixth multiplicity

\vspace{1mm}
\begin{equation}
\label{after10001qu1}
J_{T,t}^{*(i_1\ldots i_6)}={\int\limits_t^{*}}^T
\ldots
{\int\limits_t^{*}}^{t_2}
d{\bf w}_{t_1}^{(i_1)}
\ldots d{\bf w}_{t_6}^{(i_6)}
\end{equation}

\vspace{5mm}
\noindent
the following 
expansion

$$
J_{T,t}^{*(i_1\ldots i_6)}
=\hbox{\vtop{\offinterlineskip\halign{
\hfil#\hfil\cr
{\rm l.i.m.}\cr
$\stackrel{}{{}_{p\to \infty}}$\cr
}} }
\sum\limits_{j_1, \ldots, j_6=0}^{p}
C_{j_6 \ldots j_1}\zeta_{j_1}^{(i_1)}\ldots
\zeta_{j_6}^{(i_6)}
$$

\vspace{5mm}
\noindent
that converges in the mean-square sense is valid, where
$i_1, \ldots, i_6=0, 1,\ldots,m,$

\vspace{1mm}
$$
C_{j_6 \ldots j_1}=
\int\limits_t^T\phi_{j_6}(t_6)\ldots
\int\limits_t^{t_2}\phi_{j_1}(t_1)dt_1\ldots dt_6;
$$

\vspace{4mm}
\noindent
another notations are the same as in Theorems~{\rm 7--9}.}

\vspace{2mm}

Let us denote

\vspace{-1mm}
\begin{equation}
\label{str11}
I_{(l_1\ldots \hspace{0.2mm}l_k)T,t}^{*(i_1\ldots i_k)}
=
{\int\limits_t^{*}}^T
(t-t_k)^{l_k} \ldots {\int\limits_t^{*}}^{t_{2}}
(t-t_1)^{l_1} d{\bf f}_{t_1}^{(i_1)}\ldots
d{\bf f}_{t_k}^{(i_k)},
\end{equation}

\vspace{4mm}
\noindent
where $i_1,\ldots, i_k=1,\dots,m,$\ \  $l_1,\ldots,l_k=0,\ 1,\ldots.$

Below we will consider the iterated Stratonovich 
stochastic integrals (\ref{str11}) as well as the iterated
Ito stochastic integrals 
$I_{(l_1\ldots \hspace{0.2mm}l_k)T,t}^{(i_1\ldots i_k)}$
defined by (\ref{ll1}).

According to the standard relations between iterated Ito 
and Stratonovich stochastic integrals as well as according to Theorems 5, 7,
we obtain

\vspace{1mm}
$$
I_{(000)\tau_{p+1},\tau_p}^{(i_1 i_2 i_3)}=
I_{(000)\tau_{p+1},\tau_p}^{*(i_1 i_2 i_3)}
+{\bf 1}_{\{i_1=i_2\}}\frac{1}{2}I_{(1)\tau_{p+1},\tau_p}^{(i_3)}-
$$

\begin{equation}
\label{331}
-
{\bf 1}_{\{i_2=i_3\}}\frac{1}{2}\left(\Delta
I_{(0)\tau_{p+1},\tau_p}^{(i_1)}+
I_{(1)\tau_{p+1},\tau_p}^{(i_1)}\right)\ \ \ \hbox{w.\ p.\ 1},
\end{equation}

\vspace{4mm}
\noindent
where

\vspace{-2mm}
\begin{equation}
\label{dr1}
I_{(000)\tau_{p+1},\tau_p}^{*(i_1 i_2 i_3)}=
\hbox{\vtop{\offinterlineskip\halign{
\hfil#\hfil\cr
{\rm l.i.m.}\cr
$\stackrel{}{{}_{q\to \infty}}$\cr
}} }
\sum\limits_{j_1, j_2, j_3=0}^{q}
C_{j_3 j_2 j_1}\zeta_{j_1}^{(i_1)}\zeta_{j_2}^{(i_2)}\zeta_{j_3}^{(i_3)}\ \ \
(i_1, i_2, i_3=1,\ldots,m),
\end{equation}

\vspace{5mm}
\noindent
where $C_{j_3j_2j_1}$ is defined by
(\ref{hhh1}), (\ref{jjj1}).

From (\ref{331}), (\ref{dr1}) and (\ref{yyy1aaa}),
(\ref{yyy2}) we obtain 
the following approximation

\vspace{1mm}
$$
I_{(000)\tau_{p+1},\tau_p}^{(i_1i_2i_3)q}
=\sum_{j_1,j_2,j_3=0}^{q}
C_{j_3j_2j_1}
\zeta_{j_1}^{(i_1)}\zeta_{j_2}^{(i_2)}\zeta_{j_3}^{(i_3)}-
\frac{1}{4}{\bf 1}_{\{i_1=i_2\}}\Delta^{3/2}\left(
\zeta_0^{(i_3)}+\frac{1}{\sqrt{3}}\zeta_1^{(i_3)}\right)-
$$

\begin{equation}
\label{900}
-\frac{1}{4}{\bf 1}_{\{i_2=i_3\}}\Delta^{3/2}\left(
\zeta_0^{(i_1)}-\frac{1}{\sqrt{3}}\zeta_1^{(i_1)}\right).
\end{equation}

\vspace{6mm}

For the case $i_1=i_2=i_3$ it is comfortable to use the folowing 
well known relation

\begin{equation}
\label{sad003}
I_{(000)\tau_{p+1},\tau_p}^{(i_1 i_1 i_1)}
=\frac{1}{6}\Delta^{3/2}\left(
\left(\zeta_0^{(i_1)}\right)^3-3
\zeta_0^{(i_1)}\right)\ \ \ \hbox{w.\ p.\ 1}.
\end{equation}

\vspace{4mm}

Let us consider the iterated
Ito stochastic integrals

\vspace{-1mm}
$$
I_{(100)\tau_{p+1},\tau_p}^{(i_{3}i_{2}i_{1})},\ \ \
I_{(010)\tau_{p+1},\tau_p}^{(i_{3}i_{2}i_{1})},\ \ \
I_{(001)\tau_{p+1},\tau_p}^{(i_{3}i_{2}i_{1})}.
$$

\vspace{3mm}

According to the standard relations between iterated Ito 
and Stratonovich stochastic integrals as well as according to 
Theorems 5, 7, we obtain

\vspace{1mm}
$$
I_{(001)\tau_{p+1},\tau_p}^{(i_1i_2i_3)}
=I_{(001)\tau_{p+1},\tau_p}^{*(i_1i_2i_3)}+
\frac{1}{2}{\bf 1}_{\{i_1=i_2\}}
{I}_{(2)\tau_{p+1},\tau_p}^{(i_3)}+
$$

\begin{equation}
\label{d1}
+
\frac{1}{4}{\bf 1}_{\{i_2=i_3\}}\left(\Delta^2 {I}_{(0)\tau_{p+1},\tau_p}^{(i_1)}-
{I}_{(2)\tau_{p+1},\tau_p}^{(i_1)}\right)\ \ \ {\rm w.\ p.\ 1},
\end{equation}

\vspace{5mm}
$$
I_{(010)\tau_{p+1},\tau_p}^{(i_1i_2i_3)}
=
I_{(010)\tau_{p+1},\tau_p}^{*(i_1i_2i_3)}
+\frac{1}{4}{\bf 1}_{\{i_1=i_2\}}
{I}_{(2)\tau_{p+1},\tau_p}^{(i_3)}+
$$

\begin{equation}
\label{d2}
+
\frac{1}{4}{\bf 1}_{\{i_2=i_3\}}\left(\Delta^2 {I}_{(0)\tau_{p+1},\tau_p}^{(i_1}-
{I}_{(2)\tau_{p+1},\tau_p}^{(i_1)}\right)\ \ \ {\rm w.\ p.\ 1},
\end{equation}

\vspace{5mm}
$$
I_{(100)\tau_{p+1},\tau_p}^{(i_1i_2i_3)}
=I_{(100)\tau_{p+1},\tau_p}^{*(i_1i_2i_3)}
+\frac{1}{4}{\bf 1}_{\{i_1=i_2\}}
{I}_{(2)\tau_{p+1},\tau_p}^{(i_3)}-
$$

\begin{equation}
\label{d3}
-\frac{1}{2}{\bf 1}_{\{i_2=i_3\}}\left({I}_{(2)\tau_{p+1},\tau_p}^{(i_1)}+\Delta
{I}_{(1)\tau_{p+1},\tau_p}^{(i_1)}\right)\ \ \ {\rm w.\ p.\ 1},
\end{equation}

\vspace{5mm}
\noindent
where
$$
I_{(001)\tau_{p+1},\tau_p}^{*(i_1i_2i_3)}=
\hbox{\vtop{\offinterlineskip\halign{
\hfil#\hfil\cr
{\rm l.i.m.}\cr
$\stackrel{}{{}_{q\to \infty}}$\cr
}} }
\sum_{j_1,j_2,j_3=0}^{q}
C_{j_3j_2j_1}^{001}
\zeta_{j_1}^{(i_1)}\zeta_{j_2}^{(i_2)}\zeta_{j_3}^{(i_3)},
$$

\vspace{1mm}
$$
I_{(010)\tau_{p+1},\tau_p}^{*(i_1i_2i_3)}
=\hbox{\vtop{\offinterlineskip\halign{
\hfil#\hfil\cr
{\rm l.i.m.}\cr
$\stackrel{}{{}_{q\to \infty}}$\cr
}} }
\sum_{j_1,j_2,j_3=0}^{q}
C_{j_3j_2j_1}^{010}
\zeta_{j_1}^{(i_1)}\zeta_{j_2}^{(i_2)}\zeta_{j_3}^{(i_3)},
$$

\vspace{1mm}

$$
I_{(100)\tau_{p+1},\tau_p}^{*(i_1i_2i_3)}=
\hbox{\vtop{\offinterlineskip\halign{
\hfil#\hfil\cr
{\rm l.i.m.}\cr
$\stackrel{}{{}_{q\to \infty}}$\cr
}} }
\sum_{j_1,j_2,j_3=0}^{q}
C_{j_3j_2j_1}^{100}
\zeta_{j_1}^{(i_1)}\zeta_{j_2}^{(i_2)}\zeta_{j_3}^{(i_3)},
$$

\vspace{7mm}
\noindent
where
$C_{j_3j_2j_1}^{001},$
$C_{j_3j_2j_1}^{010},$
$C_{j_3j_2j_1}^{100}$
are defined by (\ref{hhh3})-(\ref{hhh5}) and
(\ref{jjj3})-(\ref{jjj5}).
From (\ref{d1})--(\ref{d3}) and (\ref{yyy1aaa}),
(\ref{yyy2}), (\ref{zzzz1})
we obtain the following approximations

\vspace{1mm}
$$
I_{(001)\tau_{p+1},\tau_p}^{(i_1i_2i_3)q}
=\sum_{j_1,j_2,j_3=0}^{q}
C_{j_3j_2j_1}^{001}
\zeta_{j_1}^{(i_1)}\zeta_{j_2}^{(i_2)}\zeta_{j_3}^{(i_3)}+
\frac{1}{6}{\bf 1}_{\{i_1=i_2\}}\Delta^{5/2}\left(
\zeta_0^{(i_3)}+\frac{\sqrt{3}}{2}\zeta_1^{(i_3)}+\frac{1}{2\sqrt{5}}\zeta_2^{(i_3)}\right)+
$$

\begin{equation}
\label{901}
+
\frac{1}{12}{\bf 1}_{\{i_2=i_3\}}\Delta^{5/2}\left(
2\zeta_0^{(i_1)}-\frac{\sqrt{3}}{2}\zeta_1^{(i_1)}-\frac{1}{2\sqrt{5}}\zeta_2^{(i_1)}\right),
\end{equation}

\vspace{5mm}
$$
I_{(010)\tau_{p+1},\tau_p}^{(i_1i_2i_3)q}
=\sum_{j_1,j_2,j_3=0}^{q}
C_{j_3j_2j_1}^{010}
\zeta_{j_1}^{(i_1)}\zeta_{j_2}^{(i_2)}\zeta_{j_3}^{(i_3)}+
\frac{1}{12}{\bf 1}_{\{i_1=i_2\}}\Delta^{5/2}\left(
\zeta_0^{(i_3)}+\frac{\sqrt{3}}{2}\zeta_1^{(i_3)}+\frac{1}{2\sqrt{5}}\zeta_2^{(i_3)}\right)+
$$

\begin{equation}
\label{902}
+
\frac{1}{12}{\bf 1}_{\{i_2=i_3\}}\Delta^{5/2}\left(
2\zeta_0^{(i_1)}-\frac{\sqrt{3}}{2}\zeta_1^{(i_1)}-\frac{1}{2\sqrt{5}}\zeta_2^{(i_1)}\right),
\end{equation}

\vspace{5mm}
$$
I_{(100)\tau_{p+1},\tau_p}^{(i_1i_2i_3)q}
=\sum_{j_1,j_2,j_3=0}^{q}
C_{j_3j_2j_1}^{100}
\zeta_{j_1}^{(i_1)}\zeta_{j_2}^{(i_2)}\zeta_{j_3}^{(i_3)}+
\frac{1}{12}{\bf 1}_{\{i_1=i_2\}}\Delta^{5/2}\left(
\zeta_0^{(i_3)}+\frac{\sqrt{3}}{2}\zeta_1^{(i_3)}+\frac{1}{2\sqrt{5}}\zeta_2^{(i_3)}\right)+
$$

\begin{equation}
\label{903}
+
\frac{1}{12}{\bf 1}_{\{i_2=i_3\}}\Delta^{5/2}\left(
\zeta_0^{(i_1)}-\frac{1}{\sqrt{5}}\zeta_2^{(i_1)}\right).
\end{equation}

\vspace{6mm}

Let us consider 
the iterated Ito stochastic integral of multiplicity 4.
According to the standard relations between iterated Ito 
and Stratonovich stochastic integrals as well as according to 
Theorems 6, 8, we get

\vspace{1mm}

$$
I_{(0000)\tau_{p+1},\tau_p}^{(i_1 i_2 i_3 i_4)}=
I_{(0000)\tau_{p+1},\tau_p}^{*(i_1 i_2 i_3 i_4)}+
\frac{1}{2}{\bf 1}_{\{i_1=i_2\}}I_{(10)\tau_{p+1},\tau_p}^{(i_3 i_4)}-
$$

\vspace{2mm}
$$
-\frac{1}{2}{\bf 1}_{\{i_2=i_3\}}\Biggl(
I_{(10)\tau_{p+1},\tau_p}^{(i_1 i_4)}-
I_{(01)\tau_{p+1},\tau_p}^{(i_1 i_4)}\Biggr)-
\frac{1}{2}{\bf 1}_{\{i_3=i_4\}}\Biggl(
\Delta I_{(00)\tau_{p+1},\tau_p}^{(i_1 i_2)}+
I_{(01)\tau_{p+1},\tau_p}^{(i_1 i_2)}\Biggr)-
$$

\vspace{2mm}
\begin{equation}
\label{uu1}
-\frac{1}{8}\Delta^2{\bf 1}_{\{i_1=i_2\}}{\bf 1}_{\{i_3=i_4\}}\ \ \ 
{\rm w.\ p.\ 1},
\end{equation}

\vspace{5mm}

$$
I_{(0000)\tau_{p+1},\tau_p}^{*(i_1 i_2 i_3 i_4)}=
\hbox{\vtop{\offinterlineskip\halign{
\hfil#\hfil\cr
{\rm l.i.m.}\cr
$\stackrel{}{{}_{q\to \infty}}$\cr
}} }
\sum\limits_{j_1, j_2, j_3, j_4=0}^{q}
C_{j_4 j_3 j_2 j_1}\zeta_{j_1}^{(i_1)}\zeta_{j_2}^{(i_2)}\zeta_{j_3}^{(i_3)}
\zeta_{j_4}^{(i_4)},
$$

\vspace{6mm}

$$
I_{(0000)\tau_{p+1},\tau_p}^{(i_1 i_2 i_3 i_4)q}=
\sum\limits_{j_1, j_2, j_3, j_4=0}^{q}
C_{j_4 j_3 j_2 j_1}\zeta_{j_1}^{(i_1)}\zeta_{j_2}^{(i_2)}\zeta_{j_3}^{(i_3)}
\zeta_{j_4}^{(i_4)}+
\frac{1}{2}{\bf 1}_{\{i_1=i_2\}}I_{(10)\tau_{p+1},\tau_p}^{(i_3 i_4)q}-
$$

\vspace{3mm}
$$
-\frac{1}{2}{\bf 1}_{\{i_2=i_3\}}\Biggl(
I_{(10)\tau_{p+1},\tau_p}^{(i_1 i_4)q}-
I_{(01)\tau_{p+1},\tau_p}^{(i_1 i_4)q}\Biggr)-
\frac{1}{2}{\bf 1}_{\{i_3=i_4\}}\Biggl(
\Delta I_{(00)\tau_{p+1},\tau_p}^{(i_1 i_2)q}+
I_{(01)\tau_{p+1},\tau_p}^{(i_1 i_2)q}\Biggr)-
$$

\vspace{3mm}
$$
-\frac{1}{8}\Delta^2{\bf 1}_{\{i_1=i_2\}}{\bf 1}_{\{i_3=i_4\}},
$$

\vspace{6mm}

\noindent
where
$$
I_{(00)\tau_{p+1},\tau_p}^{(i_1 i_2)q},\ \ \
I_{(01)\tau_{p+1},\tau_p}^{(i_1 i_2)q},\ \ \ 
I_{(10)\tau_{p+1},\tau_p}^{(i_1 i_2)q}
$$

\vspace{4mm}
\noindent
are determined by the relations  (\ref{qqqq1}), (\ref{yyy5}), (\ref{yyy6})
and $C_{j_4j_3j_2j_1}$ is defined by 
(\ref{hhh2}), (\ref{jjj2}).

For the case $i_1=i_2=i_3=i_4$ it is comfortable to use the folowing 
well known relation

\vspace{1mm}
$$
I_{(0000)\tau_{p+1},\tau_p}^{(i_1 i_1 i_1 i_1)}
=\frac{1}{24}\Delta^{2}\left(
\left(\zeta_0^{(i_1)}\right)^4-
6\left(\zeta_0^{(i_1)}\right)^2+3\right)\ \ \ \hbox{w.\ p.\ 1}.
$$

\vspace{4mm}

Let us consider 
the iterated Ito stochastic integral of fifth multiplicity 
using Theorems 6, 9

\vspace{1mm}

$$
I_{(00000)\tau_{p+1},\tau_p}^{(i_1 i_2 i_3 i_4 i_5)}=
I_{(00000)\tau_{p+1},\tau_p}^{*(i_1 i_2 i_3 i_4 i_5)}+
\frac{1}{2}{\bf 1}_{\{i_1=i_2\}}I_{(100)\tau_{p+1},\tau_p}^{(i_3 i_4 i_5)}-
$$

\vspace{2mm}
$$
-\frac{1}{2}{\bf 1}_{\{i_2=i_3\}}\Biggl(
I_{(100)\tau_{p+1},\tau_p}^{(i_1 i_4 i_5)}-
I_{(010)\tau_{p+1},\tau_p}^{(i_1 i_4 i_5)}\Biggr)-
\frac{1}{2}{\bf 1}_{\{i_3=i_4\}}\Biggl(
I_{(010)\tau_{p+1},\tau_p}^{(i_1 i_2 i_5)}-
I_{(001)\tau_{p+1},\tau_p}^{(i_1 i_2 i_5)}\Biggr)-
$$

\vspace{2mm}
$$
-\frac{1}{2}{\bf 1}_{\{i_4=i_5\}}\Biggl(
\Delta I_{(000)\tau_{p+1},\tau_p}^{(i_1 i_2 i_3)}+
I_{(001)\tau_{p+1},\tau_p}^{(i_1 i_2 i_3)}\Biggr)-
\frac{1}{8}{\bf 1}_{\{i_1=i_2\}}{\bf 1}_{\{i_3=i_4\}}
I_{(2)\tau_{p+1},\tau_p}^{(i_5)}-
$$

\vspace{2mm}
$$
-\frac{1}{8}{\bf 1}_{\{i_2=i_3\}}{\bf 1}_{\{i_4=i_5\}}\left(
\Delta^2 I_{(0)\tau_{p+1},\tau_p}^{(i_1)}+
2\Delta I_{(1)\tau_{p+1},\tau_p}^{(i_1)}+
I_{(2)\tau_{p+1},\tau_p}^{(i_1)}\right)+
$$

\vspace{2mm}
\begin{equation}
\label{uu2}
-\frac{1}{8}{\bf 1}_{\{i_1=i_2\}}{\bf 1}_{\{i_4=i_5\}}\left(
\Delta I_{(1)\tau_{p+1},\tau_p}^{(i_3)}+
I_{(2)\tau_{p+1},\tau_p}^{(i_3)}\right)\ \ \ {\rm w.\ p.\ 1},
\end{equation}

\vspace{7mm}

$$
I_{(00000)\tau_{p+1},\tau_p}^{*(i_1 i_2 i_3 i_4 i_5)}=
\hbox{\vtop{\offinterlineskip\halign{
\hfil#\hfil\cr
{\rm l.i.m.}\cr
$\stackrel{}{{}_{q\to \infty}}$\cr
}} }
\sum\limits_{j_1, j_2, j_3, j_4, j_5=0}^{q}
C_{j_5j_4 j_3 j_2 j_1}\zeta_{j_1}^{(i_1)}\zeta_{j_2}^{(i_2)}\zeta_{j_3}^{(i_3)}
\zeta_{j_4}^{(i_4)}\zeta_{j_5}^{(i_5)},
$$

\vspace{7mm}

$$
I_{(00000)\tau_{p+1},\tau_p}^{(i_1 i_2 i_3 i_4 i_5)q}=
\sum\limits_{j_1, j_2, j_3, j_4, j_5=0}^{q}
C_{j_5j_4 j_3 j_2 j_1}\zeta_{j_1}^{(i_1)}\zeta_{j_2}^{(i_2)}\zeta_{j_3}^{(i_3)}
\zeta_{j_4}^{(i_4)}\zeta_{j_5}^{(i_5)}+
\frac{1}{2}{\bf 1}_{\{i_1=i_2\}}I_{(100)\tau_{p+1},\tau_p}^{(i_3 i_4 i_5)q}-
$$

\vspace{2mm}
$$
-\frac{1}{2}{\bf 1}_{\{i_2=i_3\}}\Biggl(
I_{(100)\tau_{p+1},\tau_p}^{(i_1 i_4 i_5)q}-
I_{(010)\tau_{p+1},\tau_p}^{(i_1 i_4 i_5)q}\Biggr)-
\frac{1}{2}{\bf 1}_{\{i_3=i_4\}}\Biggl(
I_{(010)\tau_{p+1},\tau_p}^{(i_1 i_2 i_5)q}-
I_{(001)\tau_{p+1},\tau_p}^{(i_1 i_2 i_5)q}\Biggr)-
$$

\vspace{2mm}
$$
-\frac{1}{2}{\bf 1}_{\{i_4=i_5\}}\Biggl(
\Delta I_{(000)\tau_{p+1},\tau_p}^{(i_1 i_2 i_3)q}+
I_{(001)\tau_{p+1},\tau_p}^{(i_1 i_2 i_3)q}\Biggr)-
\frac{1}{8}{\bf 1}_{\{i_1=i_2\}}{\bf 1}_{\{i_3=i_4\}}I_{(2)\tau_{p+1},\tau_p}^{(i_5)}-
$$

\vspace{2mm}
$$
-\frac{1}{8}{\bf 1}_{\{i_2=i_3\}}{\bf 1}_{\{i_4=i_5\}}\left(
\Delta^2 I_{(0)\tau_{p+1},\tau_p}^{(i_1)}+
2\Delta I_{(1)\tau_{p+1},\tau_p}^{(i_1)}+
I_{(2)\tau_{p+1},\tau_p}^{(i_1)}\right)+
$$

\vspace{2mm}
$$
-\frac{1}{8}{\bf 1}_{\{i_1=i_2\}}{\bf 1}_{\{i_4=i_5\}}\left(
\Delta I_{(1)\tau_{p+1},\tau_p}^{(i_3)}+
I_{(2)\tau_{p+1},\tau_p}^{(i_3)}\right),
$$

\vspace{4mm}
\noindent
where

\vspace{-1mm}
$$
I_{(000)\tau_{p+1},\tau_p}^{(i_1 i_2 i_3)q},\ \ \
I_{(100)\tau_{p+1},\tau_p}^{(i_1 i_2 i_3)q},\ \ \
I_{(010)\tau_{p+1},\tau_p}^{(i_1 i_2 i_3)q},\ \ \
I_{(001)\tau_{p+1},\tau_p}^{(i_1 i_2 i_3)q},\ \ \
I_{(0)\tau_{p+1},\tau_p}^{(i_1)},\ \ \
I_{(1)\tau_{p+1},\tau_p}^{(i_1)},\ \ \
I_{(2)\tau_{p+1},\tau_p}^{(i_1)}
$$
 
\vspace{4mm}
\noindent
are determined by (\ref{900}), (\ref{901})--(\ref{903}),
(\ref{yyy1aaa}), (\ref{yyy2}), (\ref{zzzz1}) and
$C_{j_5j_4j_3j_2j_1}$
is defined by (\ref{hhh6}), (\ref{jjj6}).

For the case $i_1=\ldots=i_5$ it is comfortable to use the folowing 
well known relation

\vspace{1mm}
$$
I_{(00000)\tau_{p+1},\tau_p}^{(i_1 i_1 i_1 i_1 i_1)}
=\frac{1}{120}\Delta^{5/2}\left(
\left(\zeta_0^{(i_1)}\right)^5-10\left(\zeta_0^{(i_1)}\right)^3\Delta
+15\zeta_0^{(i_1)}\Delta^2\right)\ \ \ \hbox{w.\ p.\ 1}.
$$

\vspace{5mm}

Clearly, the expansions from Theorems 4--10 are simpler than the expansions 
from Theorems 1, 2.
However, the calculation of the mean-square
approximation error for the expansions from Theorems 4--10 turns out to be 
much more difficult than for the expansions
from Theorems 1, 2. We will demonstrate this fact below.

The case $k=1$ is actually not interesting. For $k=1$, 
the Ito and Stratonovich stochastic
integrals of a smooth non-random function are equal each other
w.~p.~1. 
Moreover, for $k=2$ 

$$
I_{(00)\tau_{p+1},\tau_p}^{(i_1 i_2)}=I_{(00)\tau_{p+1},\tau_p}^{*(i_1 i_2)}\ \ \ (i_1\ne i_2)\ \ \
\hbox{w.~p.~1.}
$$

\vspace{4mm}

Consider the triple Stratonovich stochastic integral defined by

\vspace{1mm}
$$
I_{(000)T,t}^{*(i_1i_2i_3)}=
\int\limits_t^{*T}\int\limits_t^{*t_{3}}
\int\limits_t^{*t_{2}}
d{\bf f}_{t_1}^{(i_1)}
d{\bf f}_{t_2}^{(i_2)}
d{\bf f}_{t_3}^{(i_3)}\ \ \ (i_1, i_2, i_3=1,\ldots,m).
$$

\vspace{4mm}                      

In view of the standard relations between 
Ito and Stratonovich stochastic integrals and
also Theorems 1, 2, 5, 7, we obtain

$$
{\sf M}\left\{\left(I_{(000)T,t}^{*(i_1i_2i_3)}-
I_{(000)T,t}^{*(i_1i_2i_3)q}\right)^2\right\}=
$$

\vspace{1mm}
$$
={\sf M}\left\{\left(I_{(000)T,t}^{(i_1i_2i_3)}+
{\bf 1}_{\{i_1=i_2\}}
\frac{1}{2}\int\limits_t^T\int\limits_t^{\tau}dsd{\bf f}_{\tau}^{(i_3)}
+{\bf 1}_{\{i_2=i_3\}}
\frac{1}{2}\int\limits_t^T\int\limits_t^{\tau}d{\bf f}_{s}^{(i_1)}d\tau-
I_{(000)T,t}^{*(i_1i_2i_3)q}\right)^2\right\}=
$$

\vspace{1mm}
$$
={\sf M}\Biggl\{\Biggl(I_{(000)T,t}^{(i_1i_2i_3)}-I_{(000)T,t}^{(i_1i_2i_3)q}+
I_{(000)T,t}^{(i_1i_2i_3)q}+
\Biggr.\Biggr.
$$

\vspace{1mm}
\begin{equation}
\label{tango3}
+
\Biggl.\Biggl.
{\bf 1}_{\{i_1=i_2\}}
\frac{1}{2}\int\limits_t^T\int\limits_t^{\tau}dsd{\bf f}_{\tau}^{(i_3)}
+{\bf 1}_{\{i_2=i_3\}}
\frac{1}{2}\int\limits_t^T\int\limits_t^{\tau}d{\bf f}_{s}^{(i_1)}d\tau-
I_{(000)T,t}^{*(i_1i_2i_3)q}\Biggr)^2\Biggr\},
\end{equation}

\vspace{3mm}

$$
I_{(000)T,t}^{(i_1i_2i_3)q}=
\sum_{j_1,j_2,j_3=0}^{q}
C_{j_3j_2j_1}\Biggl(
\zeta_{j_1}^{(i_1)}\zeta_{j_2}^{(i_2)}\zeta_{j_3}^{(i_3)}
-{\bf 1}_{\{i_1=i_2\}}
{\bf 1}_{\{j_1=j_2\}}
\zeta_{j_3}^{(i_3)}-
\Biggr.
$$
\begin{equation}
\label{tango1}
\Biggl.
-{\bf 1}_{\{i_2=i_3\}}
{\bf 1}_{\{j_2=j_3\}}
\zeta_{j_1}^{(i_1)}-
{\bf 1}_{\{i_1=i_3\}}
{\bf 1}_{\{j_1=j_3\}}
\zeta_{j_2}^{(i_2)}\Biggr),
\end{equation}

\vspace{2mm}

\begin{equation}
\label{tango2}
I_{(000)T,t}^{*(i_1i_2i_3)q}=
\sum_{j_1,j_2,j_3=0}^{q}
C_{j_3j_2j_1}
\zeta_{j_1}^{(i_1)}\zeta_{j_2}^{(i_2)}\zeta_{j_3}^{(i_3)},
\end{equation}

\vspace{5mm}
\noindent
where $I_{(000)T,t}^{(i_1i_2i_3)q}$ 
is the approximation defined by the formula (\ref{r1})
(also see (\ref{a3})) 
for the case $k=3$
and $I_{(000)T,t}^{*(i_1i_2i_3)q}$ 
is the approximation
based on Theorems 5, 7 (see (\ref{ccc7}) below).

Substituting (\ref{tango1}) and (\ref{tango2}) into (\ref{tango3}) yields

\vspace{2mm}
$$
{\sf M}\left\{\left(I_{(000)T,t}^{*(i_1i_2i_3)}-
I_{(000)T,t}^{*(i_1i_2i_3)q}\right)^2\right\}=
$$

\vspace{2mm}
$$
={\sf M}
\left\{\left(I_{(000)T,t}^{(i_1i_2i_3)}-I_{(000)T,t}^{(i_1i_2i_3)q}+
{\bf 1}_{\{i_1=i_2\}}
\left(\frac{1}{2}
\int\limits_t^T\int\limits_t^{\tau}dsd{\bf f}_{\tau}^{(i_3)}-
\sum_{j_1,j_3=0}^{q}
C_{j_3j_1j_1}
\zeta_{j_3}^{(i_3)}\right)+\right.\right.
$$

\vspace{2mm}
$$
+{\bf 1}_{\{i_2=i_3\}}\left(
\frac{1}{2}\int\limits_t^T\int\limits_t^{\tau}d{\bf f}_{s}^{(i_1)}d\tau-
\sum_{j_1,j_3=0}^{q}
C_{j_3j_3j_1}
\zeta_{j_1}^{(i_1)}\right)
\left.\left.-{\bf 1}_{\{i_1=i_3\}}
\sum_{j_1,j_2=0}^{q}
C_{j_1j_2j_1}
\zeta_{j_2}^{(i_2)}\right)^2\right\}\le
$$

\vspace{5mm}
$$
\le 4 \Biggl({\sf M}\left\{\left(I_{(000)T,t}^{(i_1i_2i_3)}-
I_{(000)T,t}^{(i_1i_2i_3)q}\right)^2\right\}+\Biggr.
{\bf 1}_{\{i_1=i_2\}}F^{(i_3)}_q+
$$

\vspace{2mm}
\begin{equation}
\label{tango4}
\Biggl.+
{\bf 1}_{\{i_2=i_3\}}G^{(i_1)}_q+
{\bf 1}_{\{i_1=i_3\}}H^{(i_2)}_q\Biggr),
\end{equation}

\vspace{5mm}
\noindent
where

$$
F^{(i_3)}_q={\sf M}\left\{\left(
\frac{1}{2}\int\limits_t^T\int\limits_t^{\tau}dsd{\bf f}_{\tau}^{(i_3)}-
\sum_{j_1,j_3=0}^{q}
C_{j_3j_1j_1}
\zeta_{j_3}^{(i_3)}\right)^2\right\},
$$

\vspace{2mm}
$$
G^{(i_1)}_q={\sf M}\left\{\left(
\frac{1}{2}\int\limits_t^T\int\limits_t^{\tau}d{\bf f}_{s}^{(i_1)}d\tau-
\sum_{j_1,j_3=0}^{q}
C_{j_3j_3j_1}
\zeta_{j_1}^{(i_1)}\right)^2\right\},
$$

\vspace{2mm}
$$
H^{(i_2)}_q={\sf M}\left\{\left(
\sum_{j_1,j_2=0}^{q}
C_{j_1j_2j_1}
\zeta_{j_2}^{(i_2)}\right)^2\right\}.
$$

\vspace{5mm}

For the cases of Legendre polynomials and 
trigonometric functions, we have the equalities
\cite{2011-2}-\cite{2018aaa}, 
\cite{2010-2}-\cite{2013}, \cite{271a},
\cite{arxiv-5}, \cite{arxiv-7}

$$
\lim\limits_{q\to\infty}F^{(i_3)}_q=0,\ \ \
\lim\limits_{q\to\infty} G^{(i_1)}_q=0,\ \ \ 
\lim\limits_{q\to\infty}H^{(i_2)}_q=0.
$$

\vspace{4mm}

However, in accordance with (\ref{tango4}) the value

$$
{\sf M}\left\{\left(I_{(000)T,t}^{*(i_1i_2i_3)}-
I_{(000)T,t}^{*(i_1i_2i_3)q}\right)^2\right\}
$$

\vspace{3mm}
\noindent
with a finite $q$ can be estimated by the sum of

\vspace{-1mm}
\begin{equation}
\label{r2}
4{\sf M}\left\{\left(I_{(000)T,t}^{(i_1i_2i_3)}-
I_{(000)T,t}^{(i_1i_2i_3)q}\right)^2\right\},
\end{equation}                                              

\vspace{3mm}
\noindent
and three additional terms of a rather complex structure. 
The value (\ref{r2}) can be calculated exactly
using Theorem 3 or estimated using (\ref{qq4}) for the case $k=3$.

As is easily observed, this peculiarity will also apply to 
the iterated Stratonovich stochastic
integrals of multiplicities 4 and 5 with the only 
difference that the number of additional terms like
$F^{(i_3)}_q$,
$G^{(i_1)}_q$, and $H^{(i_2)}_q$
will be considerably higher and their structure will be more complicated.
Therefore, the payment for relatively simple expansions of 
the iterated Stratonovich
stochastic integrals (Theorems 4--10) in comparison with the iterated
Ito stochastic integrals (Theorems 1, 2) is a much more 
difficult calculation or estimation 
procedure of their mean-square
approximation errors.

\vspace{5mm}

\section{Explicit One-Step Strong Numerical Schemes of Orders 2.0 and 2.5
Based
on the Unified Taylor--Stratonovich expansion}

\vspace{5mm}

Consider the explicit one-step strong numerical scheme of order 2.5
based on the so-called unified Taylor--Stratonovich 
expansion \cite{2006}, \cite{2017-1}-\cite{2010-1}

\vspace{2mm}
$$
{\bf y}_{p+1}={\bf y}_p+\sum_{i_{1}=1}^{m}B_{i_{1}}
\hat I_{(0)\tau_{p+1},\tau_p}^{*(i_{1})}+\Delta \bar {\bf a}
+\sum_{i_{1},i_{2}=1}^{m}G_{i_{2}}
B_{i_{1}}\hat  I_{(00)\tau_{p+1},\tau_p}^{*(i_{2}i_{1})}+
$$

\vspace{2mm}
$$
+
\sum_{i_{1}=1}^{m}\Biggl(G_{i_{1}}\bar {\bf a}\left(
\Delta \hat I_{(0)\tau_{p+1},\tau_p}^{*(i_{1})}+
\hat I_{(1)\tau_{p+1},\tau_p}^{*(i_{1})}\right)
-\bar LB_{i_{1}}\hat I_{(1)\tau_{p+1},\tau_p}^{*(i_{1})}\Biggr)+
$$

\vspace{2mm}
$$
+\sum_{i_{1},i_{2},i_{3}=1}^{m} G_{i_{3}}G_{i_{2}}
B_{i_{1}}\hat I_{(000)\tau_{p+1},\tau_p}^{*(i_{3}i_{2}i_{1})}+
\frac{\Delta^2}{2}\bar L\bar {\bf a}+
$$

\vspace{2mm}
$$
+\sum_{i_{1},i_{2}=1}^{m}
\Biggl(G_{i_{2}}\bar LB_{i_{1}}\left(
\hat I_{(10)\tau_{p+1},\tau_p}^{*(i_{2}i_{1})}-
\hat I_{(01)\tau_{p+1},\tau_p}^{*(i_{2}i_{1})}
\right)
-\bar LG_{i_{2}}B_{i_{1}}\hat I_{(10)\tau_{p+1},\tau_p}^{*(i_{2}i_{1})}
+\Biggr.
$$

\vspace{2mm}
$$
\Biggl.+G_{i_{2}}G_{i_{1}}\bar {\bf a}\left(
\hat I_{(01)\tau_{p+1},\tau_p}
^{*(i_{2}i_{1})}+\Delta \hat I_{(00)\tau_{p+1},\tau_p}^{*(i_{2}i_{1})}
\right)\Biggr)+
$$

\vspace{2mm}
$$
+
\sum_{i_{1},i_{2},i_{3},i_{4}=1}^{m}G_{i_{4}}G_{i_{3}}G_{i_{2}}
B_{i_{1}}\hat I_{(0000)\tau_{p+1},\tau_p}^{*(i_{4}i_{3}i_{2}i_{1})}+
\frac{\Delta^3}{6}LL{\bf a}+
$$

\vspace{2mm}
$$
+\sum_{i_{1}=1}^{m}\Biggl(G_{i_{1}}\bar L\bar {\bf a}\left(\frac{1}{2}
\hat I_{(2)\tau_{p+1},\tau_p}^{*(i_{1})}+
\Delta \hat I_{(1)\tau_{p+1},\tau_p}^{*(i_{1})}+
\frac{\Delta^2}{2}\hat I_{(0)\tau_{p+1},\tau_p}^{*(i_{1})}\right)\Biggr.+
$$

\vspace{2mm}
$$
\Biggl.+\frac{1}{2}\bar L 
\bar L B_{i_{1}}\hat I_{(2)\tau_{p+1},\tau_p}^{*(i_{1})}-
LG_{i_{1}}\bar {\bf a}\left(\hat I_{(2)\tau_{p+1},\tau_p}^{*(i_{1})}+
\Delta \hat I_{(1)\tau_{p+1},\tau_p}^{*(i_{1})}\right)\Biggr)+
$$

\vspace{2mm}
$$
+
\sum_{i_{1},i_{2},i_{3}=1}^m\Biggl(
G_{i_{3}}\bar LG_{i_{2}}B_{i_{1}}
\left(\hat I_{(100)\tau_{p+1},\tau_p}
^{*(i_{3}i_{2}i_{1})}-\hat I_{(010)\tau_{p+1},\tau_p}
^{*(i_{3}i_{2}i_{1})}\right)
\Biggr.+
$$

\vspace{2mm}
$$
+G_{i_{3}}G_{i_{2}}\bar LB_{i_{1}}\left(
\hat I_{(010)\tau_{p+1},\tau_p}^{*(i_{3}i_{2}i_{1})}-
\hat I_{(001)\tau_{p+1},\tau_p}^{*(i_{3}i_{2}i_{1})}\right)+
$$

\vspace{2mm}

$$
+
G_{i_{3}}G_{i_{2}}G_{i_{1}}\bar {\bf a}
\left(\Delta \hat I_{(000)\tau_{p+1},\tau_p}^{*(i_{3}i_{2}i_{1})}+
\hat I_{(001)\tau_{p+1},\tau_p}^{*(i_{3}i_{2}i_{1})}\right)
-
$$

\vspace{2mm}
$$
\Biggl.-\bar LG_{i_{3}}G_{i_{2}}B_{i_{1}}
\hat I_{(100)\tau_{p+1},\tau_p}^{*(i_{3}i_{2}i_{1})}\Biggr)+
$$

\vspace{2mm}
\begin{equation}
\label{4.470}
+\sum_{i_{1},i_{2},i_{3},i_{4},i_{5}=1}^m
G_{i_{5}}G_{i_{4}}G_{i_{3}}G_{i_{2}}B_{i_{1}}
\hat I_{(00000)\tau_{p+1},\tau_p}^{*(i_{5}i_{4}i_{3}i_{2}i_{1})},
\end{equation}

\vspace{5mm}
\noindent
where $\Delta=T/N$ $(N>1)$ is a constant (for simplicity)
step of integration,\
$\tau_p=p\Delta$ $(p=0, 1,\ldots,N)$,\
$\hat I_{(l_1\ldots l_k)s,t}^{*(i_1\ldots i_k)}$ is an 
approximation of the iterated
Stratonovich stochastic integral (\ref{str11}),

\vspace{1mm}
$$
\bar{\bf a}({\bf x},t)={\bf a}({\bf x},t)-
\frac{1}{2}\sum\limits_{j=1}^m G_jB_j({\bf x},t),
$$

\vspace{2mm}
$$
\bar L=L-\frac{1}{2}\sum\limits_{j=1}^m G_0^{(j)}G_0^{(j)}=
\frac{\partial }{\partial t}+
\sum\limits_{j=1}^n \bar {\bf a}^{(j)}({\bf x},t)
\frac{\partial }{\partial {\bf x}^{(j)}},
$$

\vspace{2mm}
$$
L= {\partial \over \partial t}
+ \sum^ {n} _ {i=1} {\bf a}_i ({\bf x},  t) 
{\partial  \over  \partial  {\bf  x}_i}
+ {1\over 2} \sum^ {m} _ {j=1} \sum^ {n} _ {l,i=1}
B_{lj} ({\bf x}, t) B_{ij} ({\bf x}, t) {\partial
^{2} \over \partial {\bf x}_l \partial {\bf x}_i},
$$

\vspace{2mm}
$$
G_i = \sum^ {n} _ {j=1} B_{ji} ({\bf x}, t)
{\partial  \over \partial {\bf x}_j},\ \ \ i=1,\ldots,m,
$$

\vspace{4mm}
\noindent
$l_1,\ldots, l_k=0, 1, 2,$\ \
$i_1,\ldots, i_k=1,\ldots,m,$\ \ $k=1, 2,\ldots, 5$,\ \
$B_i$ and $B_{ij}$ are the $i$th column and the $ij$th
component of the matrix function $B$,
${\bf a}_i$ is the $i$th component of the vector function ${\bf a},$
${\bf x}_i$ is the $i$th component
of the column ${\bf x}$, 
the functions  

\vspace{-1mm}
$$
B_{i_{1}},\  \bar {\bf a},\ G_{i_{2}}B_{i_{1}},\
G_{i_{1}}\bar {\bf a},\ \bar LB_{i_{1}},\ G_{i_{3}}G_{i_{2}}B_{i_{1}},\ 
\bar L\bar {\bf a},\ LL{\bf a},\
G_{i_{2}}\bar LB_{i_{1}},\ 
$$

\vspace{-3mm}
$$
\bar LG_{i_{2}}B_{i_{1}},\ G_{i_{2}}G_{i_{1}}\bar{\bf a},\
G_{i_{4}}G_{i_{3}}G_{i_{2}}B_{i_{1}},\ G_{i_{1}}\bar L\bar{\bf a},\
\bar L\bar LB_{i_{1}},\ \bar LG_{i_{1}}\bar {\bf a},\ 
G_{i_{3}}\bar LG_{i_{2}}B_{i_{1}},\
G_{i_{3}}G_{i_{2}}\bar LB_{i_{1}},\
$$

\vspace{-3mm}
$$
G_{i_{3}}G_{i_{2}}G_{i_{1}}\bar {\bf a},\
\bar LG_{i_{3}}G_{i_{2}}B_{i_{1}},\ 
G_{i_{5}}G_{i_{4}}G_{i_{3}}G_{i_{2}}B_{i_{1}}
$$

\vspace{3mm}
\noindent
are calculated at the point $({\bf y}_p,p).$

Under the standard conditions \cite{KlPl2}, \cite{2006} the numerical 
scheme (\ref{4.470}) has 
strong order 2.5 of convergence. 
The major emphasis below will be placed on the 
approximation of the iterated
Stratonovich stochastic integrals appearing in (\ref{4.470}). 
Therefore, among 
the mentioned standard conditions, we note the
appro\-xi\-ma\-ti\-on condi\-ti\-on for 
these stochastic integrals \cite{KlPl2}, \cite{2006}, 
which has the form

\vspace{-1mm}
\begin{equation}
\label{ors}
{\sf M}\left\{\biggl(I_{(l_{1}\ldots\hspace{0.2mm} l_{k})\tau_{p+1},\tau_p}
^{*(i_{1}\ldots i_{k})} 
- \hat I_{(l_{1}\ldots\hspace{0.2mm} l_{k})\tau_{p+1},
\tau_p}^{*(i_{1}\ldots i_{k})}
\biggr)^2\right\}\le C\Delta^{6},
\end{equation}

\vspace{3mm}
\noindent
where constant $C$ is independent of
$\Delta$.

Note that if we exclude from (\ref{4.470}) the terms starting from the
term $\Delta^3 LL{\bf a}/6$, then we will have the explicit 
one-step strong numerical scheme of order 2.0 \cite{KlPl2}, 
\cite{2006}, \cite{2017-1}-\cite{2010-1}.

Using the numerical scheme (\ref{4.470}) or its modifications based 
on the Taylor--Stratonovich expansion \cite{KlPl1},
the implicit or multistep analogues of (\ref{4.470}) can be constructed
\cite{KlPl2}, \cite{2006}, \cite{2017-1}-\cite{2010-1}. The set of the
iterated Stratonovich 
stochastic integrals to be approximated for implementing 
these modifications is the same
as for the numerical scheme (\ref{4.470}) itself.
Interestingly, the truncated unified Taylor--Stratonovich expansion (the 
foundation of the numerical
scheme (\ref{4.470})) contains 12 different types of the 
iterated Stratonovich
stochastic integrals 
(\ref{str11}), which cannot be
interconnected by linear relations \cite{2006}, \cite{2017-1}-\cite{2010-1}. 
The analogous 
Taylor--Stratonovich expansion \cite{KlPl2}, \cite{KlPl1} contains
17 different types of iterated Stratonovich 
stochastic integrals, part of which 
are interconnected by linear relations
and part of which have a higher multiplicity than the iterated 
Stratonovich stochastic integrals (\ref{str11}). This
fact well explains the use of the numerical scheme (\ref{4.470}).

One of the main problems arising in the implementation of the 
numerical scheme (\ref{4.470}) is the joint
numerical modeling of the iterated Stratonovich stochastic integrals 
figuring in (\ref{4.470}). Let us consider an efficient 
numerical modeling method for 
the iterated Stratonovich stochastic integrals based on Theorems 4--10.

Using Theorems 4--9 and multiple Fourier--Legendre series,
we obtain the following 
approximations of the iterated Stratonovich stochastic 
integrals from (\ref{4.470}) \cite{2006}-\cite{arxiv-6}

\vspace{2mm}
\begin{equation}
\label{ccc1}
I_{(0)\tau_{p+1},\tau_p}^{*(i_1)}=\sqrt{\Delta}\zeta_0^{(i_1)},
\end{equation}

\vspace{2mm}

\begin{equation}
\label{ccc2}
I_{(1)\tau_{p+1},\tau_p}^{*(i_1)}=
-\frac{{\Delta}^{3/2}}{2}\left(\zeta_0^{(i_1)}+
\frac{1}{\sqrt{3}}\zeta_1^{(i_1)}\right),
\end{equation}

\vspace{2mm}

\begin{equation}
\label{ccc3}
{I}_{(2)\tau_{p+1},\tau_p}^{*(i_1)}=
\frac{\Delta^{5/2}}{3}\left(
\zeta_0^{(i_1)}+\frac{\sqrt{3}}{2}\zeta_1^{(i_1)}+
\frac{1}{2\sqrt{5}}\zeta_2^{(i_1)}\right),
\end{equation}

\vspace{3mm}
\begin{equation}
\label{ccc4}
I_{(00)\tau_{p+1},\tau_p}^{*(i_1 i_2)q}=
\frac{\Delta}{2}\left(\zeta_0^{(i_1)}\zeta_0^{(i_2)}+\sum_{i=1}^{q}
\frac{1}{\sqrt{4i^2-1}}\left(
\zeta_{i-1}^{(i_1)}\zeta_{i}^{(i_2)}-
\zeta_i^{(i_1)}\zeta_{i-1}^{(i_2)}\right)\right),
\end{equation}

\vspace{5mm}

$$
I_{(01)\tau_{p+1},\tau_p}^{*(i_1 i_2)q}=
-\frac{\Delta}{2}
I_{(00)\tau_{p+1},\tau_p}^{*(i_1 i_2)q}
-\frac{{\Delta}^2}{4}\Biggl(
\frac{1}{\sqrt{3}}\zeta_0^{(i_1)}\zeta_1^{(i_2)}+\Biggr.
$$

\vspace{2mm}
\begin{equation}
\label{ccc5}
+\Biggl.\sum_{i=0}^{q}\Biggl(
\frac{(i+2)\zeta_i^{(i_1)}\zeta_{i+2}^{(i_2)}
-(i+1)\zeta_{i+2}^{(i_1)}\zeta_{i}^{(i_2)}}
{\sqrt{(2i+1)(2i+5)}(2i+3)}-
\frac{\zeta_i^{(i_1)}\zeta_{i}^{(i_2)}}{(2i-1)(2i+3)}\Biggr)\Biggr),
\end{equation}

\vspace{5mm}

$$
I_{(10)\tau_{p+1},\tau_p}^{*(i_1 i_2)q}=
-\frac{\Delta}{2}I_{(00)\tau_{p+1},\tau_p}^{*(i_1 i_2)q}
-\frac{\Delta^2}{4}\Biggl(
\frac{1}{\sqrt{3}}\zeta_0^{(i_2)}\zeta_1^{(i_1)}+\Biggr.
$$

\vspace{2mm}
\begin{equation}
\label{ccc6}
+\Biggl.\sum_{i=0}^{q}\Biggl(
\frac{(i+1)\zeta_{i+2}^{(i_2)}\zeta_{i}^{(i_1)}
-(i+2)\zeta_{i}^{(i_2)}\zeta_{i+2}^{(i_1)}}
{\sqrt{(2i+1)(2i+5)}(2i+3)}+
\frac{\zeta_i^{(i_1)}\zeta_{i}^{(i_2)}}{(2i-1)(2i+3)}\Biggr)\Biggr),
\end{equation}

\vspace{5mm}

\begin{equation}
\label{ccc7}
I_{(000)\tau_{p+1},\tau_p}^{*(i_1 i_2 i_3)q}
=\sum_{j_1, j_2, j_3=0}^{q}
C_{j_3 j_2 j_1}
\zeta_{j_1}^{(i_1)}\zeta_{j_2}^{(i_2)}\zeta_{j_3}^{(i_3)},
\end{equation}

\vspace{3mm}

\begin{equation}
\label{ccc8}
I_{(100)\tau_{p+1},\tau_p}^{*(i_1 i_2 i_3)q}
=\sum_{j_1, j_2, j_3=0}^{q}
C_{j_3 j_2 j_1}^{100}
\zeta_{j_1}^{(i_1)}\zeta_{j_2}^{(i_2)}\zeta_{j_3}^{(i_3)},
\end{equation}

\vspace{3mm}

\begin{equation}
\label{ccc9}
I_{(010)\tau_{p+1},\tau_p}^{*(i_1 i_2 i_3)q}
=\sum_{j_1, j_2, j_3=0}^{q}
C_{j_3 j_2 j_1}^{010}
\zeta_{j_1}^{(i_1)}\zeta_{j_2}^{(i_2)}\zeta_{j_3}^{(i_3)},
\end{equation}

\vspace{3mm}

\begin{equation}
\label{ccc10}
I_{(001)\tau_{p+1},\tau_p}^{*(i_1 i_2 i_3)q}
=\sum_{j_1, j_2, j_3=0}^{q}
C_{j_3 j_2 j_1}^{001}
\zeta_{j_1}^{(i_1)}\zeta_{j_2}^{(i_2)}\zeta_{j_3}^{(i_3)},
\end{equation}

\vspace{3mm}

\begin{equation}
\label{ccc11}
I_{(0000)\tau_{p+1},\tau_p}^{*(i_1 i_2 i_3 i_4)q}
=\sum_{j_1, j_2, j_3, j_4=0}^{q}
C_{j_4 j_3 j_2 j_1}
\zeta_{j_1}^{(i_1)}\zeta_{j_2}^{(i_2)}\zeta_{j_3}^{(i_3)}
\zeta_{j_4}^{(i_4)},
\end{equation}

\vspace{3mm}

\begin{equation}
\label{ccc12}
I_{(00000)\tau_{p+1},\tau_p}^{*(i_1 i_2 i_3 i_4 i_5)q}=
\sum\limits_{j_1, j_2, j_3, j_4, j_5=0}^{q}
C_{j_5j_4 j_3 j_2 j_1}\zeta_{j_1}^{(i_1)}\zeta_{j_2}^{(i_2)}\zeta_{j_3}^{(i_3)}
\zeta_{j_4}^{(i_4)}\zeta_{j_5}^{(i_5)},
\end{equation}

\vspace{6mm}

\noindent
where the Fourier--Legendre coefficients

\vspace{1mm}
$$
C_{j_3 j_2 j_1},\ \ \ C_{j_3 j_2 j_1}^{100},\ \ \
C_{j_3 j_2 j_1}^{010},\ \ \
C_{j_3 j_2 j_1}^{001},\ \ \ 
C_{j_4 j_3 j_2 j_1},\ \ \
C_{j_5j_4 j_3 j_2 j_1}
$$ 

\vspace{4mm}
\noindent
are determined by  
(\ref{hhh1})--(\ref{hhh6}),
(\ref{jjj1})--(\ref{jjj6}).

On the basis of 
the presented 
expansions (see (\ref{ccc1})--(\ref{ccc12})) of 
iterated Stratonovich stochastic integrals we 
can see that increasing of multiplicities of these integrals 
or degree indexes of their weight functions 
leads
to increasing 
of smallness orders with respect to $\Delta$ in the mean-square sense 
for iterated stochastic integrals. This leads to a sharp decrease  
of member 
quantities (the numbers $q$)
in expansions of iterated Stratonovich stochastic 
integrals, which are required for achieving the acceptable accuracy
of approximation.
Generally speaking, the minimum values $q$ that guarantee the fulfillment
of the condition 
(\ref{ors})
for each approximation (\ref{ccc1})--(\ref{ccc12})
are different and abruptly decreasing with the growth of 
smallness order (with respect to $\Delta$) of
the approximations 
of iterated stochastic integrals.

From Theorem 3 for the case $i_1\ne i_2$ we obtain

\vspace{1mm}
$$
{\sf M}\left\{\left(I_{(00)\tau_{p+1},\tau_p}^{*(i_1 i_2)}-
I_{(00)\tau_{p+1},\tau_p}^{*(i_1 i_2)q}
\right)^2\right\}=\frac{\Delta^2}{2}
\sum\limits_{i=q+1}^{\infty}\frac{1}{4i^2-1}\le 
$$

\vspace{2mm}
\begin{equation}
\label{teac}
\le \frac{\Delta^2}{2}\int\limits_{q}^{\infty}
\frac{1}{4x^2-1}dx
=-\frac{\Delta^2}{8}{\rm ln}\left|
1-\frac{2}{2q+1}\right|\le C_1\frac{\Delta^2}{q},
\end{equation}

\vspace{5mm}
\noindent
where constant $C_1$ does not depend on $\Delta$.

As was mentioned above,
the value $\Delta$ plays the role of integration step 
in the numerical procedures for Ito stochastic differential equations. 
Thus this value is a sufficiently small.
Keeping in mind this circumstance, it is easy to notice that there 
exists such a constant $C_2$ that

\vspace{1mm}
\begin{equation}
\label{teac3}
{\sf M}\left\{\left(I_{(l_1\ldots l_k)\tau_{p+1},\tau_p}^{*(i_1\ldots i_k)}-
I_{(l_1\ldots l_k)\tau_{p+1},\tau_p}^{*(i_1\ldots i_k)q}\right)^2\right\}
\le C_2 {\sf M}\left\{\left(I_{(00)\tau_{p+1},\tau_p}^{*(i_1 i_2)}-
I_{(00)\tau_{p+1},\tau_p}^{*(i_1 i_2)q}\right)^2\right\},
\end{equation}

\vspace{4mm}
\noindent
where $I_{(l_1\ldots l_k)\tau_{p+1},\tau_p}^{*(i_1\ldots i_k)q}$
is an approximation of the iterated Stratonovich stochastic integral 
(\ref{str11}).

From (\ref{teac}) and (\ref{teac3}) we finally have

\vspace{1mm}
\begin{equation}
\label{teac4}
{\sf M}\left\{\left(I_{(l_1\ldots l_k)\tau_{p+1},\tau_p}^{*(i_1\ldots i_k)}-
I_{(l_1\ldots l_k)\tau_{p+1},\tau_p}^{*(i_1\ldots i_k)q}\right)^2\right\}
\le K \frac{\Delta^2}{q},
\end{equation}

\vspace{4mm}
\noindent
where constant $K$ does not depend on $\Delta.$ 

The same idea can be found in \cite{KlPl2} 
for the case of trigonometric functions.
Note that, in contrast to the estimate (\ref{teac4}), 
the constant $C$ in Theorems 7--9 does not depend on $q.$

We can get significantly more information about numbers $q$ 
using a different approach.
Applying the standard relation between iterated Ito
and Stratonovich  stochastic integrals, we have

\vspace{-1mm}
$$
I_{(l_1\ldots l_k)\tau_{p+1},\tau_p}^{*(i_1\ldots i_k)}=
I_{(l_1\ldots l_k)\tau_{p+1},\tau_p}^{(i_1\ldots i_k)}\ \ \ \hbox{w.~p.~1}
$$

\vspace{5mm}
\noindent
for pairwise different $i_1,\ldots,i_k=1,\ldots,m$.

Then for $i_1\ne i_2$ the following mean-square errors

\vspace{2mm}
$$
{\sf M}\biggl\{\left(I_{(00)\tau_{p+1},\tau_p}^{*(i_1 i_2)}-
I_{(00)\tau_{p+1},\tau_p}^{*(i_1 i_2)q}
\right)^2\biggr\},\ \ \
{\sf M}\biggl\{\left(I_{(10)\tau_{p+1},\tau_p}^{*(i_1 i_2)}-
I_{(10)\tau_{p+1},\tau_p}^{*(i_1 i_2)q}
\right)^2\biggr\},
$$

\vspace{3mm}
$$
{\sf M}\biggl\{\left(I_{(01)\tau_{p+1},\tau_p}^{*(i_1 i_2)}-
I_{(01)\tau_{p+1},\tau_p}^{*(i_1 i_2)q}
\right)^2\biggr\}
$$

\vspace{6mm}
\noindent
are defined by (\ref{xxx1}), (\ref{xxx2}).

Moreover, for pairwise different $i_1,\ldots,i_5=1,\ldots,m$ 
from (\ref{qq1}) we obtain

\vspace{2mm}

$$
{\sf M}\left\{\left(
I_{(01)\tau_{p+1},\tau_p}^{*(i_1i_2)}-
I_{(01)\tau_{p+1},\tau_p}^{*(i_1i_2)q}\right)^2\right\}=
\frac{\Delta^{4}}{4}-\sum_{j_1,j_2=0}^{q}
\left(C_{j_2j_1}^{01}\right)^2,
$$

\vspace{4mm}
$$
{\sf M}\left\{\left(
I_{(10)\tau_{p+1},\tau_p}^{*(i_1i_2)}-
I_{(10)\tau_{p+1},\tau_p}^{*(i_1i_2)q}\right)^2\right\}=
\frac{\Delta^{4}}{12}-\sum_{j_1,j_2=0}^{q}
\left(C_{j_2j_1}^{10}\right)^2,
$$

\vspace{4mm}

$$
{\sf M}\left\{\left(
I_{(000)\tau_{p+1},\tau_p}^{*(i_1i_2 i_3)}-
I_{(000)\tau_{p+1},\tau_p}^{*(i_1i_2 i_3)q}\right)^2\right\}=
\frac{\Delta^{3}}{6}-\sum_{j_3,j_2,j_1=0}^{q}
C_{j_3j_2j_1}^2,
$$

\vspace{4mm}
$$
{\sf M}\left\{\left(
I_{(0000)\tau_{p+1},\tau_p}^{*(i_1i_2 i_3 i_4)}-
I_{(0000)\tau_{p+1},\tau_p}^{*(i_1i_2 i_3 i_4)q}\right)^2\right\}=
\frac{\Delta^{4}}{24}-\sum_{j_1,j_2,j_3,j_4=0}^{q}
C_{j_4j_3j_2j_1}^2,
$$

\vspace{4mm}
$$
{\sf M}\left\{\left(
I_{(100)\tau_{p+1},\tau_p}^{*(i_1i_2 i_3)}-
I_{(100)\tau_{p+1},\tau_p}^{*(i_1i_2 i_3)q}\right)^2\right\}=
\frac{\Delta^{5}}{60}-\sum_{j_1,j_2,j_3=0}^{q}
\left(C_{j_3j_2j_1}^{100}\right)^2,
$$

\vspace{4mm}
$$
{\sf M}\left\{\left(
I_{(010)\tau_{p+1},\tau_p}^{*(i_1i_2 i_3)}-
I_{(010)\tau_{p+1},\tau_p}^{*(i_1i_2 i_3)q}\right)^2\right\}=
\frac{\Delta^{5}}{20}-\sum_{j_1,j_2,j_3=0}^{q}
\left(C_{j_3j_2j_1}^{010}\right)^2,
$$

\vspace{4mm}
$$
{\sf M}\left\{\left(
I_{(001)\tau_{p+1},\tau_p}^{*(i_1i_2 i_3)}-
I_{(001)\tau_{p+1},\tau_p}^{*(i_1i_2 i_3)q}\right)^2\right\}=
\frac{\Delta^5}{10}-\sum_{j_1,j_2,j_3=0}^{q}
\left(C_{j_3j_2j_1}^{001}\right)^2,
$$

\vspace{4mm}
$$
{\sf M}\left\{\left(
I_{(00000)\tau_{p+1},\tau_p}^{*(i_1 i_2 i_3 i_4 i_5)}-
I_{(00000)\tau_{p+1},\tau_p}^{*(i_1 i_2 i_3 i_4 i_5)q}\right)^2\right\}=
\frac{\Delta^{5}}{120}-\sum_{j_1,j_2,j_3,j_4,j_5=0}^{q}
C_{j_5 i_4 i_3 i_2 j_1}^2.
$$

\vspace{7mm}

For example \cite{2006}-\cite{2013},

\vspace{3mm}

$$
{\sf M}\left\{\left(
I_{(000)\tau_{p+1},\tau_p}^{*(i_1i_2 i_3)}-
I_{(000)\tau_{p+1},\tau_p}^{*(i_1i_2 i_3)6}\right)^2\right\}=
\frac{\Delta^{3}}{6}-\sum_{j_3,j_2,j_1=0}^{6}
C_{j_3j_2j_1}^2
\approx
0.01956000\Delta^3,
$$

\vspace{4mm}
$$
{\sf M}\left\{\left(
I_{(100)\tau_{p+1},\tau_p}^{*(i_1i_2 i_3)}-
I_{(100)\tau_{p+1},\tau_p}^{*(i_1i_2 i_3)2}\right)^2\right\}=
\frac{\Delta^{5}}{60}-\sum_{j_1,j_2,j_3=0}^{2}
\left(C_{j_3j_2j_1}^{100}\right)^2
\approx
0.00815429\Delta^5,
$$

\vspace{4mm}
$$
{\sf M}\left\{\left(
I_{(010)\tau_{p+1},\tau_p}^{*(i_1i_2 i_3)}-
I_{(010)\tau_{p+1},\tau_p}^{*(i_1i_2 i_3)2}\right)^2\right\}=
\frac{\Delta^{5}}{20}-\sum_{j_1,j_2,j_3=0}^{2}
\left(C_{j_3j_2j_1}^{010}\right)^2
\approx
0.01739030\Delta^5,
$$

\vspace{4mm}
$$
{\sf M}\left\{\left(
I_{(001)\tau_{p+1},\tau_p}^{*(i_1i_2 i_3)}-
I_{(001)\tau_{p+1},\tau_p}^{*(i_1i_2 i_3)2}\right)^2\right\}=
\frac{\Delta^5}{10}-\sum_{j_1,j_2,j_3=0}^{2}
\left(C_{j_3j_2j_1}^{001}\right)^2
\approx 0.02528010\Delta^5,
$$

\vspace{4mm}
$$
{\sf M}\left\{\left(
I_{(0000)\tau_{p+1},\tau_p}^{*(i_1i_2i_3 i_4)}-
I_{(0000)\tau_{p+1},\tau_p}^{*(i_1i_2i_3 i_4)2}\right)^2\right\}=
\frac{\Delta^{4}}{24}-\sum_{j_1,j_2,j_3,j_4=0}^{2}
C_{j_4 j_3 j_2 j_1}^2\approx
0.02360840\Delta^4,
$$

\vspace{4mm}
$$
{\sf M}\left\{\left(
I_{(00000)\tau_{p+1},\tau_p}^{*(i_1i_2i_3i_4 i_5)}-
I_{(00000)\tau_{p+1},\tau_p}^{*(i_1i_2i_3i_4 i_5)1}\right)^2\right\}=
\frac{\Delta^5}{120}-\sum_{j_1,j_2,j_3,j_4,j_5=0}^{1}
C_{j_5 i_4 i_3 i_2 j_1}^2\approx
0.00759105\Delta^5.
$$

\end{document}